\documentclass{amsart}
\usepackage{epsf,amscd,amssymb}
\usepackage{graphicx}

\newtheorem{thm}{Theorem}
\newtheorem{lem}{Lemma}
\newtheorem{cor}{Corollary}
\newtheorem{prop}{Proposition}

\theoremstyle{definition}

\theoremstyle{remark}
\newtheorem{rem}{Remark}
\newtheorem{exmp}{Example}

\newcommand\bR{{\bf R}}

\newcommand\bZ{{\bf Z}}

\newcommand\Hom{{\rm Hom}}

\newcommand\SI{{\bf S}}

\newcommand\clo{{\rm Cl}}

\newcommand\ra{\rightarrow}

\newcommand\emp{\emptyset}
\newcommand\eps{\epsilon}

\newcommand\vth{\vartheta}
\newcommand\vpi{\varphi}

\newcommand\sgn{{\rm sign}}
\newcommand\id{{\rm id}}

\begin{document}


\title{Geometric structures on orbifolds and
holonomy representations}
\author{Suhyoung Choi}
\thanks{The author gratefully acknowledges partial support 
by the grant number 1999-2-101-002-3 
from the Interdisciplinary Research Program of KOSEF
} 
\email{shchoi@math.snu.ac.kr}
\address{
Department of Mathematics,
College of Natural Sciences, 
Seoul National University, 
151--742 Seoul, Korea}
\date{January 2, 2003} 
\subjclass{Primary 57M50; Secondary 53A20, 53C15}
\keywords{orbifold,
deformation space, geometric structure, group representation}

\begin{abstract}
An orbifold is a topological space modeled on quotient 
spaces of a finite group actions. We can define 
the universal cover of an orbifold and the fundamental
group as the deck transformation group. 
Let $G$ be a Lie group acting on a space $X$. 
We show that the space of isotopy-equivalence classes 
of $(G,X)$-structures on a compact orbifold $\Sigma$ is 
locally homeomorphic to the space of representations of 
the orbifold fundamental group of $\Sigma$ to $G$
following the work of Thurston, Morgan, and Lok. This implies 
that the deformation space of $(G, X)$-structures on
$\Sigma$ is locally homeomorphic to 
the space of representations of 
the orbifold fundamental group to $G$ when restricted to 
the region of proper conjugation action by $G$.
\end{abstract}  

\maketitle 


\section{Introduction}

An orbifold is a topological space
with neighborhoods modeled on the orbit-spaces of
finite group actions on open subsets of euclidean spaces. 
Often orbifolds arise as quotient spaces of manifolds by proper 
actions of discrete groups.
They are so-called good orbifolds.
The manifold itself could be chosen to be a simply connected one. 
In this case, the manifold is said to be the {\em universal covering 
space} of the orbifold and the deck transformation group
the {\em {\rm (} orbifold {\rm )} fundamental group} of the orbifold. 

For example, the quotient spaces of
the hyperbolic spaces by discrete subgroups of
isometries are orbifolds (especially, when
the group has torsion elements).
Very good orbifolds are orbifolds which are quotient spaces of 
manifolds by finite group actions. The good orbifolds
are as good as manifolds since they admit universal
covering space as manifolds. 

Looking at quotients of manifolds by 
group actions as orbifolds sometimes gives us useful 
methods such as decomposition or putting geometric 
structures by cut-and-paste methods. This is one 
of the reasons why we study orbifolds instead of
just manifolds. 

Let $G$ be a Lie group acting on a space $X$ transitively
and effectively.Then $(G, X)$ is said to be a {\em geometry}.
A {\em $(G, X)$-structure} on an orbifold 
$M$ is given by a maximal atlas of charts to orbit spaces 
of finite subgroups of $G$ acting on open subsets of $X$.
A $(G,X)$-structure on an orbifold implies that 
the orbifold is good as first observed by Thurston. 
This generalizes the 
notion of $(G,X)$-structures on manifolds introduced 
by Ehresmann. 

We will define in this paper, the space 
$\mathcal{S}(\Sigma)$ of isotopy-equivalence
classes of $(G,X)$-structures on a given orbifold $\Sigma$.
Given a $(G, X)$-structure on $\Sigma$, we 
can define an immersion $D$ from its universal cover 
to $X$ and a homomorphism $h:\pi_1(\Sigma) \ra G$
for the orbifold fundamental group $\pi_1(\Sigma)$ of $\Sigma$.
$D$ is said to be a {\em developing map}, and 
$h$ a {\em holonomy homomorphism}.  
$(D, h)$ determines the $(G,X)$-structure 
but given a $(G, X)$-structure $(D, h)$ is 
determined up to the following action
\[ (D, h(\cdot)) \mapsto 
(\vth \circ D, \vth \circ h(\cdot) \circ \vth^{-1})\]
for $\vth \in G$. 
The so-called development pair
$(D, h)$ is essentially defined by an analytic continuation
of charts as in the manifold cases. (See Goldman \cite{Grep:88}
for more details on $(G, X)$-structures on manifolds.)
The space $\mathcal{S}(\Sigma)$ can be considered 
as the space of equivalence classes 
of development pairs of $(G, X)$-structures on $\tilde \Sigma$
under the isotopy action of $\tilde \Sigma$ commuting with
the deck-transformation group. 
The deformation space $\mathcal{D}(\Sigma)$
of $(G, X)$-structures on $\Sigma$
is obtained from $\mathcal{S}(\Sigma)$ 
as a quotient space by the above action of $G$. 

We need to assume that $\Sigma$ is compact and 
$\pi_1(\Sigma)$ is finitely-presented. 
One can define a map, so-called pre-holonomy map, 
\[\mathcal{PH}:\mathcal{S}(\Sigma) \ra \Hom(\pi_1(\Sigma), G)\]
induced by an isotopy-invariant function 
assigning a $(G,X)$-structure with a developing map $D$
to its holonomy homomorphism associated with $D$. 
$\Hom(\pi_1(\Sigma), G)$ is naturally a real algebraic subset 
of $G^n$ where $n$ is the number of generators of $\pi_1(\Sigma)$ 
defined by relations and hence is a topological space. 
\begin{thm} 
Let $(X, G)$ be a geometry, and
$\Sigma$ be a compact orbifold with a finitely-presented 
orbifold fundamental group $\pi_1(\Sigma)$. Then
\[\mathcal{PH}: \mathcal{S}(\Sigma) \ra 
\Hom(\pi_1(\Sigma), G)\]
is a local homeomorphism. 
\end{thm}
If the orbifold is given an additional ``cellular structure''
in some sense, then the compactness should imply 
that the fundamental group is finitely generated. 
(However, we do not prove it here.) 

The proof of the above theorem for {\em manifolds} was first given
by Thurston (perhaps much earlier by Ehresmann), again by 
Canary-Epstein-Green \cite{CER:87}, and simultaneously by 
Lok (following Morgan).
Our proof generalizes that 
in the manifold case by Lok \cite{Lok:84} 
following Morgan's lectures (see Weil \cite{We:60} and 
Canary-Epstein-Green \cite{CER:87} also). 
(We mention that this 
can be also done using Goldman's idea in \cite{Grep:88}.) 
(There are related works by Kapovich \cite{Kapo:00}
and Gallo-Kapovich-Marden \cite{GKM:00}
where some results were proved for $2$-orbifolds partially.)

Let us denote by $\Hom(\pi_1(\Sigma), G)^{s}$ the 
part where the conjugation action of $G$ given
by 
\[ h(\cdot) \mapsto g h(\cdot) g^{-1}, g \in G\] 
is stable (see \cite{Johnson}). 
Let $\mathcal{D}^{s}(\Sigma)$ 
the inverse image of the above set by $\mathcal{PH}$. 
(We assume here that $G$ is a group of $\bR$-points of 
an algebraic group defined over the real number fields.)

\begin{cor} 
Let $\Sigma$ be a compact orbifold with a finitely-presented 
orbifold fundamental group $\pi_1(\Sigma)$.
Then the map 
\[ \mathcal{H}: \mathcal{D}^s(\Sigma)/G \ra \Hom(\pi_1(\Sigma), G)^{s}/G\]
induced by $\mathcal{PH}| \mathcal{D}^s(\Sigma)$
is a local homeomorphism to its image.
\end{cor}

This result will be used in Choi-Goldman \cite{CG:00}, 
which is the main reason why we wrote this paper. 

This paper intentionally is technical and 
gives many details since 
excellent intuitive writings on the subjects are already 
available in Chapter 5 of Thurston's notes \cite{Thnote} and 
the paper by Scott \cite{Scott:83}. (See also Kato \cite{Kato:87} 
and Matsumoto and Montesinos-Amilibia \cite{MatMont}.)
Although this material has already been exposed well, 
we believe that it should be a glad duty of 
mathematical community to continuously reinterpret 
old ideas and make precise and refine
what initial attempts to understand have left to the next 
generation. Hopefully, this writing will convey
the idea of orbifolds to readers in a more rigorous way.

In fact some of the preliminary materials in this paper 
have been also exposed in  Haefliger \cite{Hae:90},
Bridson-Haefliger \cite{Hae:99},
and Ratcliffe \cite{Rat:94} 
but they study only orbifolds with geometric 
structures where $G$ acts as an isometry group on 
a space $X$. For our purposes
with geometric structures without metrics, 
their results are not enough and 
hence the need for writing them down arose. 
For geometric structures such as projectively flat and 
conformally flat structures on $n$-dimensional 
orbifolds, the results in this paper
are new except for the work of Kapovich described above.

In Section two, we introduce orbifolds, 
orbifold-maps, isotopies of orbifold-maps,
and so on. We also explain the Riemannian metric on
orbifolds and coverings by normal neighborhoods. 

In Section three, we first 
review fiber-products of topological covering spaces. 
We discuss covering orbifolds of orbifolds, and discuss its 
simple properties. We define fiber-products of orbifold-covering maps
by first doing so for orbifold-coverings of elementary neighborhoods, 
and then extending the definitions to any collection of orbifold-coverings.
We prove a theorem of Thurston that there exists a so-called 
universal covering orbifold for any orbifold by fiber-products.
We provide a proof of the fact that the deck transformation
group acts transitively on the universal covering orbifolds.  
From these results, we obtain most properties of orbifold-coverings
similar to topological coverings. Finally, we show that a good
orbifold has a manifold as the universal covering orbifold. 
The author faithfully follows and gives selective details of 
Chapter 5 of Thurston \cite{Thnote}.
(In fact, this material should 
be published by Thurston in his next book. This part has 
grown unintentionally large and is only there to provide 
a technical background to this paper.)

In Section four, we discuss the geometric structures
on orbifolds. We show that orbifolds with geometric 
structures are good, and find the developing maps
and the holonomy homomorphisms for orbifolds.
We define the deformation space of $(G,X)$-structures
on an orbifold, which is the space of equivalence classes of 
$(G,X)$-structures under isotopies and $(G,X)$-diffeomorphisms.
The so-called isotopy-equivalence 
space of $(G,X)$-structures on an orbifold $\Sigma$ 
is defined to be the space of equivalence classes of
a pair $(D, f)$ where $D$ is a developing map
for a $(G,X)$-orbifold $M$, and $f$ is a lift of 
an orbifold-diffeomorphism defined on the universal cover of 
$\Sigma$. The equivalence relation is given by 
an isotopy action on $f$. We show that $G$ acts on 
the isotopy-equivalence space so that the quotient space is 
the deformation space of $(G,X)$-structures here. 
We define a pre-holonomy map from the isotopy-equivalence
space of $(G,X)$-structures on an orbifold $M$ to the space of
representations $\Hom(\pi_1(\Sigma), G)$ given 
by sending $(D, f)$ to the holonomy homomorphism
composed with the isomorphism $\pi_1(\Sigma) \ra \pi_1(M)$ 
induced by $f$. Here $\pi_1(\Sigma)$ denotes the 
deck transformation group of $\Sigma$.

In Section five, we prove Theorem 1 that the space of 
isotopy classes of $(G,X)$-structures on an orbifold is locally
homeomorphic to the space of representations 
of the fundamental group to $G$ by the pre-holonomy map. 
The proof is essentially the same as that of 
Morgan and Lok \cite{Lok:84} but we modify 
slightly for clarity and completeness. The basic idea 
is to deform small neighborhoods first
and patch them together using ``bump'' functions
as we change the representation by a small amount
in a cone-neighborhood of the representation space
as described in Canary-Epstein-Green \cite{CER:87}.

The local finite group actions complicate the proof
somewhat but not greatly if we use the old ideas of
Palais-Stewart \cite{PS:60} which yield three 
necessary lemmas on conjugating 
finite group actions by diffeomorphisms 
in the beginning of Section five. For technical purposes, 
we explain the Riemannian metric structures on orbifolds
here. The proof of Theorem 1  is as follows: 
\begin{itemize}
\item[(I)] We choose three model-neighborhood
coverings $\{U_i\}, \{W_i\},$ and $\{V_i\}$ 
of the orbifold so that they are nested, i.e., 
\[\clo(U_i) \subset V_i, \clo(V_i) \subset W_i 
\hbox{ for each } i\]  
Then we lift each open set to a connected open 
set in the universal cover and choose deck 
transformations so that patching the lifted
sets by a selection of 
deck transformations gives us back 
the orbifold. 
\item[(II)] We show that there is
a local section of the pre-holonomy map:
we show that as we deform 
the holonomy representations, 
we can deform the model neighborhoods 
by conjugating with respect to finite group action
deformations. We patch the deformations together
to form a deformation of $M$. 
\item[(III)] We finally show that the pre-holonomy
map is a local homeomorphism; that is, 
if the holonomy homomorphisms of two $(G, X)$-structures 
are equal and they are close with one another, then 
we will find a $(G,X)$-orbifold-diffeomorphism between 
them.
\end{itemize}

We would like to thank Yves Benoist, William Goldman, Karsten Grove, 
Silvio Levi, Misha Kapovich, Hyuk Kim, Inkang Kim, 
John Millson, and Shmuel Weinberger 
for their helpful comments and encouragements. 
We especially thank Bill Thurston for discovering the ideas in 
this paper. We thank greatly the referee for his 
detailed comments to polish this paper up. 

\section{Topology of orbifolds}

In this paper, we assume that the action of 
a group on a topological space is locally faithful; that is, 
for each nonidentity element $g$ restricts to
nonidentity on each open subset of the 
space where the group acts on. 
By this requirement, the set of fixed 
points of any nontrivial subgroup is always nowhere dense.
Also, if two elements agree locally, then 
they are equal. For finite groups, this is
always true by M.H.A. Newman \cite{New:31}.

In this paper, we will consider only differentiable maps and 
sets with differentiable structures. Most of the difficulty of
topological group actions disappear in this case. 
For example, Newman's result is trivial in differentiable cases. 

There are extensive literature on the finite group actions
on manifolds using many interesting methods 
for which this author has no expertise on: 
For example, many actions on an $n$-cell are not conjugate to 
linear actions. In fact there are finite group actions without fixed points
(see Floyd and Richardson \cite{FloydR} and
Buchdhal, Kwasik, and Schultz \cite{BKS}).
However, if one chooses a sufficiently small ball around 
a fixed point, then the action of the finite group fixing 
that point is conjugate to a linear action.
Also, every smooth action of a compact Lie group on 
a $3$-dimensional Euclidean space is differentiably 
conjugate to a linear action as claimed by Thurston 
and shown by Kwasik and Schultz \cite{KwasikS}. 
(See Davis \cite{Davis} for a survey.)

Most of the material of this section 
is in Chapter 5 of Thurston (91/12/19 version) \cite{Thnote}:
We will repeat it here for reader's convenience 
and the difficulty of the writing there and 
some omissions. See also Satake \cite{Satake:56}, \cite{Satake:57}, 
Kato \cite{Kato:87}, and 
Matsumoto and Montesinos-Amilibia \cite{MatMont}).

An $n$-dimensional {\em orbifold} is a Hausdorff, second-countable space
$Y$ so that each point has a neighborhood homeomorphic to
the quotient space $U$ of an open set in $\bR^n$ by an action of
a finite group. Moreover, if such a neighborhood
$V$ of $y$, modeled on a pair $(\tilde V, G_1)$
is a subset of another such neighborhood
$U$, modeled on a pair $(\tilde U, G_2),$ 
then the inclusion map $\phi_{V, U}: V \ra U$ lifts to 
an imbedding $\tilde \phi_{V, U}: \tilde V \ra \tilde U$
equivariant with respect to a homomorphism
$\psi_{V, U}: G_1 \ra G_2$ so that the following diagram
is commutative.
\begin{eqnarray}\label{eqn:orbdefn}
\tilde V &\stackrel{\tilde \phi_{V, U}}{\longrightarrow} &\tilde U
\nonumber \\ 
\downarrow & & \downarrow \nonumber \\
\tilde V/G_1 &\longrightarrow & \tilde U/\psi_{V, U}(G_1) \nonumber \\  
\downarrow & & 
\begin{array}{c} \downarrow\\ \tilde U/G_2 \\ \downarrow \end{array}
\nonumber \\
V &\stackrel{\phi_{V, U}}{\longrightarrow} & \quad U 
\end{eqnarray} 
Note that the pair $(\tilde \phi_{V, U}, \psi_{V, U})$ 
can be chosen differently;
i.e., the pair $\vth \circ \tilde \phi_{V, U}$ and 
$\vth\circ \psi_{V, U}(\cdot)\circ \vth^{-1}$ for
$\vth \in G_2$ satisfies the above equation as well.
Thus, $(\tilde \phi_{V, U}, \psi_{V, U})$
associated to the map $\phi_{V, U}$ is unique up to an element of $G_2$. 
If $\phi_{V, U}:\tilde V \hookrightarrow \tilde U$ and 
$\phi_{U, W}:\tilde U\hookrightarrow \tilde W$ are inclusion maps, then
we are forced to have
\begin{eqnarray}
\tilde \phi_{V, W} & = &\vth \circ \tilde \phi_{U, W} \circ
\tilde \phi_{V, U} \hbox{ and } \nonumber \\
\psi_{V, W}(\cdot) &=& \vth \circ \psi_{U, W}\circ \psi_{V, U}(\cdot)
\circ \vth^{-1}, \mbox{ for } \vth \in G_3
\end{eqnarray}
where $G_3$ is the finite group associated with $W$. 
($V$ is said to be a {\em model neighborhood}
and $(\tilde V, G_V)$ the {\em model pair} where 
$G_V$ is a finite group acting on an open subset 
$\tilde V$ of $\bR^n$.)

A maximal family of coverings $\mathcal{O}$ with models satisfying
the above conditions is said to be an {\em orbifold structure}
on $Y$. (That is, an orbifold structure is 
a maximal collection of model pairs with 
inclusion equivalence classes 
satisfying the above properties.)
$Y$ is said to be the {\em underlying space}
of $(Y, \mathcal{O})$ where $\mathcal{O}$ has models as above.
Given an orbifold $M$, we denote by $X_M$ the underlying space 
in this paper. (We won't often distinguish between 
the underlying space and the orbifold itself, particularly later
on.)

Given two orbifolds $M$ and $N$, an {\em orbifold-map}
is a map $f:X_M \ra X_N$ so that for each point $x$ 
of $X_N$, a neighborhood of $x$ modeled on
$(U, G)$, and an inverse image of $y$, there is 
a neighborhood of $y$ modeled on $(V, G')$ and
a smooth map $\tilde f:V \ra U$ inducing $f$ equivariant
with respect to a homomorphism $\psi:G' \ra G$.
(That is, we record the lifting $\tilde f$ but $\tilde f$ 
is determined only up to $G$, $G'$, i.e., 
the map $g\circ \tilde f \circ g'', g \in G, g'' \in G'$
and homomorphism $\psi(\cdot)$ changed 
to $g\circ \psi(g''(\cdot)g^{\prime \prime -1}) \circ g^{-1}$.
Moreover, such liftings have to be consistent 
in a way that one can take two copies
of equation (\ref{eqn:orbdefn}) for $M$ and $N$ 
and write $\tilde f$ and induced maps between 
corresponding elements.)

An {\em orbifold-diffeomorphism} is an orbifold-map 
which has an inverse function with lifts
that again forms an orbifold-map.

An {\em orbifold with boundary} is a Hausdorff,
second-countable space so that each point 
has a neighborhood modeled on an open set
intersected with the upper-half space and a finite group
acting on it. The {\em interior} is a set of
points with neighborhoods modeled on open balls.
The {\em boundary} is the complement of the interior.
(The boundary is a boundaryless orbifold of codimension one.)

A {\em singular} point $x$ of an orbifold is a point
of the underlying space which has a neighborhood
whose model neighborhood has a nontrivial element of 
the group fixing a point corresponding to $x$. 
A nonsingular point, so-called regular point, of 
an orbifold always has a neighborhood homeomorphic to
a ball. The set of regular points is an open dense subset
of the underlying space since the set of fixed points 
of a finite group in a model pair is a nowhere dense 
closed set. The set of singular points is nowhere dense
since so is the set of fixed points of a differentiable 
group action. In this paper, often a point in the open subset of 
the model pair is said to be {\em regular} or {\em singular} 
depending on what the image in the quotient is.

A {\em suborbifold} of an orbifold $N$ is an imbedded subset
$Y$ of $X_N$ with an orbifold structure so that for each point 
$x$ of $Y$, and a neighborhood $V$ modeled on $(V', G)$,
the neighborhood $V \cap Y$ is modeled on 
$(V'\cap P, G|P)$ where $P$ is a submanifold of
$\bR^n$ where $G$ acts, and $G|P$ denotes the image 
subgroup of the restriction homomorphism 
to groups acting on $P$. 

The boundary of an orbifold is a suborbifold clearly. 

\begin{exmp}
A class of examples are given as follows:
Let $M$ be a manifold and $\Gamma$ a discrete group
acting on $M$ properly but not necessarily freely.
Then $M/\Gamma$ has an orbifold structure:
Let $x$ be a point of $M/\Gamma$ and 
$\tilde x$ a point of $M$ corresponding to $x$.
Then a subgroup $I_{\tilde x}$ of $\Gamma$ fixes $\tilde x$.
There is a ball-neighborhood $U$ of $\tilde x$ 
on which $I_{\tilde x}$ acts and for any 
$g \in \Gamma - I_{\tilde x}$, $g(U) \cap U$ is empty. 
Then $U/I_{\tilde x}$ is a neighborhood of $x$
modeled on $(U, I_{\tilde x})$. If $V$ is another
such neighborhood in $U/ I_{\tilde x}$ containing a point $y$,
then a component $V'$ of its inverse image in $U$
is acted upon by a subgroup $I'$ of $I_{\tilde x}$. 
Also, for any $g\in \Gamma - I'$, $g(V')\cap V'$
is empty. Therefore, the inclusion $V \ra U/I_{\tilde x}$ 
satisfies the conditions for equations (\ref{eqn:orbdefn}).
\end{exmp}

Given two orbifolds $M$ and $N$, the product space
$X_M \times X_N$ obviously has an orbifold structure;
i.e., we model on $(U\times V, G_U \times G_V)$
if $(U, G_U)$ and $(V, G_V)$ are model pairs 
for neighborhoods of $M$ and $N$ respectively. 
The product space with this orbifold structure is
denote by $M \times N$.

A {\em homotopy} of two orbifold-maps 
$f_1, f_2: M\ra N$ from an orbifold $M$ to another one $N$
is an orbifold-map $F: M\times [0,1] \ra N$ where 
$[0,1]$ is the unit interval and $F(x, 0) = f_1(x)$ and
$F(x, 1) = f_2(x)$ for every $x\in M$.
We define an orbifold-map $F_t:M \ra N$ to be given by
$F_t(x) = F(x, t)$ with appropriate liftings in
model pairs of $M$ and $N$. 

Given an orbifold $M$, an {\em isotopy} $f:M \ra M$ is
a self-orbifold-diffeomorphism, i.e., an automorphism, so that there is 
a homotopy $F:M\times [0,1] \ra M$ so that 
$F_0$ is the identity map and $F_1 = f$,
and $F_t$ is an orbifold-diffeomorphism for each $t$.

Two orbifold-diffeomorphisms $f_1, f_2:M \ra M'$
are {\em isotopic} if there is a homotopy
$F:M \times [0,1] \ra M'$ so that 
$F_0 = f_1$ and $F_1 = f_2$ and 
$F_t$ are orbifold-diffeomorphisms. 

A {\em Riemannian} metric on an orbifold 
is a Riemannian metric on each model open set
invariant under the associated finite group action
and where inclusion induced maps for model pairs 
are isometries. 
(See Satake \cite{Satake:56} and \cite{Satake:57}
for more details.)
A {\em partition of unity} for an open cover of 
an orbifold is a collection of functions whose supports 
are in compact subsets of the elements of the open cover, 
which sum to $1$ and correspond to finite-group-invariant smooth functions 
on model pairs.
An orbifold with an open cover has a partition of unity.
We can always put a Riemannian metric on a compact orbifold:
Cover the orbifold by the modeled neighborhoods 
and choose a locally finite subcover $\{V_i\}$ and a partition of unity.
Let $(U_i, G_i)$ be the model pairs.
Choose a Riemannian metric on $U_i$ 
and by taking an average over the finite
group action $G_i$, we obtain an invariant metric on
each modeled neighborhood $V_i$. 
Let $\{\phi_i\}$ be the partition of unity on the orbifold
such that the support of $\phi_i$ is in a compact subset of $V_i$
for each $i$. Then $\phi_i$ pulls back to a smooth function 
$\tilde \phi_i$ which is $G_i$ equivariant.
Let us choose a Riemannian metric $\mu_i$ on 
each $U_i$ so that $G_i$ acts by isometries. 
Then on each $U_i$, we may form a smooth pseudo metric  
$\tilde \phi_j \mu_j^*$ which is induced from the inclusion map to $U_j$.
Clearly, $\mu'_i = \sum_j \tilde \phi_j \mu_j^*$ is a smooth metric on $U_i$
where $G_i$ acts as isometries 
and the inclusion maps induce isometries.
This defines a global metric on the orbifold. 
Furthermore, for any model neighborhood, its model pair 
has a well-defined induced Riemannian metric 
invariant under the group action.
We may consider such metrics Riemannian metrics on 
orbifolds 

We make a quotient space of the tangent bundle $T(U_i)$ over
$U_i$ by $G_i$ to obtain $2n$-dimensional orbifold
$O_i$. We can easily patch $O_i$s together to 
obtain a $2n$-orbifold $T(M)$ with a map $p:T(M) \ra M$ 
so that the inverse image of a point is 
a vector space modulo a finite group action.
Let $T_{x_0}(M)$ denote the fiber over $x_0 \in M$.

If $x_0 \in M$ is a singular point in $V_i$, then we can choose
an open ball $U_{x_0}$ in $U_i$ so that the subgroup 
$G_{x_0}$ of $G_i$
fixing the point $\tilde x_0$ 
corresponding to $x_0$ acts on it. 
Then there is a neighborhood
of $V_{x_0}$ of $x_0$ which is modeled on 
$(U_{x_0}, G_{x_0})$.

An exponential map from $T_{x_0}(M)$ to $V_{x_0}$ is locally defined 
by the exponential map on the model open set $U_{x_0}$ which
is clearly invariant under the finite group action
if $x_0$ is singular. If $x_0$ is regular, we can
use the ordinary exponential map. We can obviously 
patch these maps to obtain a global map
$\exp_{x_0}: T_{x_0}(M) \ra M$. 

We can find $r>0$ so that $\exp_{\tilde x_0}$ 
imbeds the ball $B_r(0) \subset T_{\tilde x_0}U_{x_0}$ of
positive radius $ < r$ to a strictly convex ball in $U_{x_0}$.
(They have smooth convex boundary.)
Thus, the exponential map from each $x_0 \in M$
sends a quotient space of a ball positive radius $< r$
to a quotient space of a strictly
convex ball in $M$. The images are said to be 
{\em normal neighborhoods}.

\begin{lem} 
Let $x$ be a point of an orbifold $M$. 
Then $x$ has a model pair $(\tilde V, G_V)$ so 
that $\tilde V$ is simply-connected. 
\end{lem}
\begin{proof} 
We may choose $\tilde V$ sufficiently small so that 
$G_V$ action is conjugate to a linear action.
For example, we can choose the normal neighborhood.
This is Proposition 5.4.1 of \cite{Thnote}. 
\end{proof}

We will need to find a covering of an orbifold by model sets of
certain forms. A covering $\{O_i\}$ of an orbifold $M$
is said to be a {\em nice covering} if it satisfies 
\begin{itemize}
\item Each $O_i$ is connected and open.
\item $O_i$ has a model pair $(\tilde O_i, G_i)$ so that 
$\tilde O_i$ is simply-connected.
\item The intersection of any finite collection of $O_i$ 
has the above two properties.
\end{itemize} 

\begin{prop}\label{prop:nicecover}
An orbifold $M$ has a nice locally finite covering.
\end{prop}
\begin{proof}
We will assume that $M$ is compact. If $M$ is only locally compact,
the proof is similar.
First cover $M$ by a finite collection of normal neighborhoods. 
For a point $x$ of a model pair $(U, G_U)$ of a normal neighborhood, 
there exists a radius $r> 0$, such that the ball $B_r(x)$ of 
radius $r$ in $U$ has the property that for each pair of 
points $y$ and $z$ of $B_r(x)$, there exists a unique geodesic
segment connecting $y$ and $z$. Thus, using the Lebesgue number,
we can find a real number $r_0 >0$ so that for each point $x$ of $M$, 
the ball $B_r(x)$ of radius $r$, $0 < r < r_0$ in $M$ 
has the generic convexity property;
i.e., for each generic pair of points $y$ and $z$ of $B_r(x)$, 
there exists a unique geodesic connecting $y$ and $z$.
Given two balls $B_{r/8}(x)$ and $B_{r/8}(y)$, their intersection 
$B_{r/8}(x) \cap B_{r/8}(y)$, which is open, can have only
one component. By an induction, we can show that 
any finite intersection of balls $\bigcap_{i=1}^n B_{r/8}(x_i)$ 
is connected. The collection consisting of $B_{r/8}(x_i)$, $i=1, \dots, n$,
covering $M$ is a nice covering.
\end{proof}

\section{Fiber products and the universal covering orbifolds}

In some cases, we will allow covering spaces to have many 
components for convenience, which does not cause too much confusion.
In this case, a {\em morphism} of 
topological covering spaces $(X_1, p_1)$ and $(X_2, p_2)$
is a map $f: X_1 \ra X_2$ such that $p_1 = p_2 \circ f$ so that
$f$ induces injective correspondence between the components of
$X_1$ and components of $X_2$ and if $X_1$ is connected, and 
$X_1$ and $X_2$ have base points, we require that $f$ sends $X_1$ to 
the component of $X_2$ containing the base point. Note that  
$f$ need not be surjective. However, if $X_1$ and $X_2$ are connected,
morphisms are surjective.

Let us briefly review the notion of (topological) fiber-products in 
ordinary covering space 
theory so that we can generalize the notion to that for 
orbifold-covering spaces:
Given a sequence of covering maps $p_i:Y_i \ra Y$
for $i$ in some index set $I$, 
in the ordinary sense, one can form a fiber-product 
$p^f:Y^f\ra Y$ by setting 
$Y^f$ to be the subspace of $\prod_{i\in I} Y_i$ such that 
for $(x_i)_{i \in I} \in Y^f$
\[p_j\circ \pi_j((x_i)_{i\in I}) = p_k\circ \pi_k((x_i)_{i\in I})\]
for all $j, k$ and $\pi_i:\prod_{i\in I} Y_i \ra Y_i$ the projection 
to the $i$-th factor. The covering map $p^f:Y^f \ra Y$
is given by $p^f((x_i)) = p_1(x_1)$,
and $Y^f$ covers $Y_i$ by a morphism $p'_i:Y^f \ra Y_i$, 
so that $p_i\circ p'_i = p^f$, which is given by 
the projection to the $i$th-factor.

Given base-points $y_i$ of $Y_i$ mapping to a base point $y$ of $Y$, 
we decide that the corresponding point $y^f = (y_i)_{i\in I}$ 
of $Y^f$ be a base point of $Y^f$.

The fiber product is not necessarily connected. 
For example, consider the fiber-product of two-fold covering
of and a four-fold covering of $\SI^1$.
The fiber-product maps to $\SI^1$ by $8$ to $1$, 
but it has two components mapping $4$ to $1$. 

Let $Y^f_0$ be the component containing the base point. 
Then $Y^f_0$ is the covering of $Y$ corresponding to 
$\bigcap_{i \in I} p_{i, *}(\pi_1(Y_i, y_i))$. 
This follows by considering which loops of 
$Y$ based at $y$ can be lifted to loops in the fiber-product;
A loop lifts if and only if it is in the above intersection.
We also see that each of the other components 
is a covering corresponding to 
\[\bigcap_{i \in I} p_{i, *}(\pi_1(Y_i, y'_i))\]
where $(y'_i)_{i \in I}$ is a point of the component
so that $p_i(y'_i) = y$. 

We can also verify that 
$Y^f_0$ has a {\em universal property} that 
if $(Y'', p'')$ is a connected covering space of 
$Y$ with a base point $y''$ 
and $q_i:Y'' \ra Y_i$ is a covering morphism
sending $y''$ to $y_i$ for each $i$,
then there exists a covering morphism
$q':Y'' \ra Y^f_0$ sending $y''$ to $y^f$.
(Of course, we can change $y_i$ to any $y'_i$ 
with $p_i(y'_i) = y$.)

Also, the universal property characterizes $Y^f_0$
up to covering isomorphisms. 
That is, if $(Z, p_Z)$ is a connected covering space of $Y$ so that 
there is a covering morphism $q_{Z, i}:Z \ra Y_i$,
so that $z \mapsto y_i$ for a base point $z$ of $Z$ and each $i$,
and $Z$ satisfies the universal property of $Y^f_0$ above, 
then there exists a covering isomorphism 
$L: Z \ra Y^f_0$ such that 
\begin{equation} q_{Z, i} = p'_i\circ L. \end{equation}

Finally, if the collection $\{(Y_i, p_i)\}$ contains all 
the covering spaces of $Y$ up to isomorphisms, 
we see that components of $Y^f$ are universal covers of $Y$.

\begin{exmp}\label{exmp:fiberprod}
Since we need the group action descriptions to define 
orbifold-fiber products, 
we view above example more in terms of group actions:
Let a connected manifold $Y$ have a connected regular covering space 
$\tilde Y$ with the covering map $\tilde p$
and subgroups $\Gamma_i$s of the deck transformation 
group $\Gamma$.
Then let $Y_i$ be the quotient space $\tilde Y/\Gamma_i$
and $p_i:Y_i \ra Y$ the covering map for each $i$.
Let $\Gamma_i \backslash \Gamma$ denote the right coset
space of $\Gamma_i$ in $\Gamma$.
The projection map
\[\tilde p_i:\tilde Y \times \Gamma_i\backslash \Gamma
\ra \tilde Y\]
induces a map
\[(\tilde Y \times \Gamma_i\backslash \Gamma)/\Gamma
\ra \tilde Y / \Gamma\]
where $\Gamma$ acts on the first space
by \[\gamma(\tilde x, \Gamma_i \gamma_i)
= (\gamma(\tilde x), \Gamma_i \gamma_i \gamma^{-1})
\hbox{ for } \gamma \in \Gamma\]
and on the second space in the usual way. 
Since $\Gamma$ acts transitively on each sheets of 
\[\tilde Y \times \Gamma_i \backslash \Gamma\]
and $\Gamma_i$ acts on the sheet 
$\tilde Y \times \Gamma_i \cdot 1$, we see that
this map is obviously the same map as $p_i$. 
The fiber product of $\tilde p_i$s is 
clearly equal to the projection
\[\tilde Y \times \prod_{i\in I} \Gamma_i\backslash \Gamma 
\ra \tilde Y.\]
Define the left action of $\Gamma$ on the first space by 
\[\gamma(\tilde x, (\Gamma_i \gamma_i)_{i\in I}) = 
(\gamma(\tilde x), (\Gamma_i\gamma_i\gamma^{-1})_{i\in I})
\hbox{ for } \gamma \in \Gamma.\]

Since $\Gamma$-equivalence classes of the first space 
correspond exactly to the fiber products of 
the $\Gamma$-equivalence classes of $\tilde Y$, 
the fiber product of $p_i$s equals 
\[(\tilde Y \times \prod_{i\in I} (\Gamma_i\backslash \Gamma))/
\Gamma \ra \tilde Y /\Gamma \]
induced by the projection.
As before, the fiber product may have many components.
Using path-considerations again, we see that 
each component is isomorphic to a cover 
\[\tilde Y/\bigcap_{i \in I} \gamma_i\Gamma_i \gamma_i^{-1}\]
of $\tilde Y/\Gamma$ where $\gamma_i$, $i \in I$, 
is a sequence of a coset-representative of 
$\Gamma_i \backslash \Gamma$ for each $i$.
Another way to see this is that given a component 
\[\tilde Y \times \prod_{i\in I} \Gamma_i \gamma_i\]
for a sequence $\gamma_i$, the group acting on 
it equals $\bigcap_{i \in I} \gamma_i \Gamma_i \gamma_i^{-1}$. 
Thus, the quotient of the component 
corresponds to the component described above. 

If $\{\Gamma_i\}_{i \in I}$ are all of the subgroups of $\Gamma$, 
then each component equals $\tilde Y$. 

We note that the covering map from the fiber-product 
to $Y$ is given by 
sending $[\tilde x, (\Gamma_i \gamma_i)_{i\in I}]$ to 
$[\tilde p(\tilde x)]$,
and the covering morphism from the fiber product to 
$Y_i$ for each $i$ is given by
sending $[\tilde x, (\Gamma_i\gamma_i)_{i\in I}]$ 
to $[\gamma_i(\tilde x)]$ in $\tilde Y/\Gamma_i$. 
\end{exmp}

A {\em covering orbifold}\/ of an orbifold $M$ is an orbifold
$\tilde M$ with a surjective orbifold-map $p:X_{\tilde M} \ra X_M$ such that 
each point $x \in X_M$ has a neighborhood $U$,
so-called an elementary neighborhood, with
a homeomorphism $\phi:\tilde U/G_U \ra U$ and 
an open subset of $\tilde U$ in $\bR^n$ or $\bR^{n, +}$ with
a group $G_U$ acting on it,  so that 
each component $V_i$ of $p^{-1}(U)$ has a homeomorphism
$\tilde \phi_i:  \tilde U/G_i \ra V_i$ (in the orbifold
structure) where $G_i$ is a subgroup of $G_U$.
We require the quotient map $\tilde U \ra V_i$ induced by 
$\tilde \phi_i$ composed with $p$ should be the quotient map 
$\tilde U \ra U$ induced by $\phi$.
(We often don't assume $X_{\tilde M}$ is connected.
In this case, only components of $M$ need to be orbifolds
and $M$ itself does not.
See also Bridson-Haefliger \cite{Hae:99})

A {\em fiber} of a point of $M$ is the inverse image $p^{-1}(x)$.

Given an orbifold-map $f:Y \ra Z$ and a covering $(Z_1, p_1)$  
of $Z$, if an orbifold-map $\tilde f:Y \ra Z_1$ satisfies
$p_1\circ \tilde f = f$ and $\tilde f$ lifts for
every model pair of points of $Z_1$ in the consistent way for $Z$,
$\tilde f$ is said to be a {\em lifting} of $f$. 

Two covering orbifolds $(Y_1, p_1)$ and $(Y_2, p_2)$ of 
an orbifold $Y$ are {\em isomorphic} if
there is an orbifold-diffeomorphism $f:Y_1 \ra Y_2$
so that $p_2 \circ f = p_1$.
A {\em covering automorphism} or {\em deck transformation}
$Y_1 \ra Y_1$ is a covering isomorphism of $Y_1$ to itself. 
(Thus, $f$ is a lifting of $p_1$.)
More generally, a {\em morphism} $(Y_1, p_1) \ra (Y_2, p_2)$
is an orbifold-map $f:Y_1\ra Y_2$ so that $p_2\circ f= p_1$
where distinct components go to distinct components
(see Proposition \ref{prop:automorphisms}).  
If $Y_1$ and $Y_2$ are connected, then $f$ obviously is surjective
by an open and closedness argument.
For orbifolds with base points, we require that 
a morphism should preserve base points.
A covering $(Y_1, p_1)$ is {\em regular} if the automorphism 
group acts transitively on fibers over regular points. 
Given coverings $(Y_1, p_1)$ over $Y$ and 
$(Y_2, p_2)$ over $Z$, an orbifold-map $f:Y_1 \ra Y_2$ 
{\em covers} an orbifold-map $g:Y \ra Z$ if the following 
diagram is commutative:
\begin{eqnarray*}
Y_1 &\stackrel{f}{\longrightarrow} & Y_2 \\
p_1 \downarrow & & \downarrow p_2 \\
Y &\stackrel{g}{\longrightarrow} & Z. 
\end{eqnarray*}

In this paper, if $Z$ is a cover of an orbifold $Y$, 
then by $Z^r$ we mean the inverse image of the regular
part $Y^r$ of $Y$. The inverse image consists of 
regular points of $Z$; however, the converse is 
not true in general. (It will be clear from the context 
whether one means just a regular part or the part over
the regular part.)

We now go over some preliminary results on 
orbifold-covering maps:
\begin{prop}\label{prop:nonsingular}
Let $(Y_1, p_1)$ and $(Y_2, p_2)$ be coverings     
over an orbifold $Y$. 
Let $f:Y_1 \ra Y_2$ be a covering morphism
so that $f:Y_1^r \ra Y_2^r$ is a covering isomorphism
where $Y_1^r$ and $Y_2^r$ are inverse images of
the regular part $Y^r$ of $Y$.
Then $f$ itself is a covering isomorphism.
\end{prop}
\begin{proof}
For the model pairs, the groups have to be isomorphic.
The rest is straightforward. 
\end{proof}

\begin{prop}\label{prop:samemap} 
Let $(Y_1, p_1)$, $(Y_2, p_2)$, and $(Y_3, p_3)$ 
be coverings over $Y$. Let $f_1:Y_1 \ra Y_3$, 
$f_2: Y_2 \ra Y_3$, and $f_3: Y_1 \ra Y_2$ 
be covering morphisms 
so that $f_1| Y_1^r = f_2\circ f_3| Y_1^r$. 
Then $f_1 = f_2\circ f_3$. 
\end{prop}
\begin{proof} 
Again a local consideration proves this
in a straightforward manner.
\end{proof}

By a {\em path}, we mean a smooth orbifold-map from an interval to 
an orbifold.
The two central properties of covering space theory 
survive in the orbifold-covering space theory: 
\begin{prop}
Let $Y$ be an orbifold and $p:Y' \ra Y$ an orbifold-covering map.
\begin{itemize}
\item Let $x$ be a point of $Y$ and $x'$ a point in $p^{-1}(x)$.
A path $f: I \ra Y$ such that $f(0) = x$ 
lifts to a unique path $f':I \ra Y'$ in $Y'$ 
such that $f'(0) = x'$. 
\item Let 
$f_1:Z \ra Y'$ and $f_2: Z \ra Y'$ be orbifold-maps 
lifting $f:Z \ra Y$. If $f_1(x) = f_2(x)$ for a regular point $x \in Z$, 
then $f_1 = f_2$. 
\end{itemize}
\end{prop}
\begin{proof} 
The first one has the same proof as the ordinary topological 
theory since liftings are determined up to the action of 
local groups and we only need to match them. 
We use the open and closeness for the second one.
\end{proof}

We discuss somewhat about so-called doubling. 
A {\em mirror point} is a singular point of an orbifold 
with a model pair $(U, G)$ where $G$ is an order-two group
acting on the open subset $U$ of $\bR^n$ fixing a hyperplane 
meeting $U$. The set of mirror points is said to 
be the {\em mirror set}.

We can form an orbifold $M^d$ covering $M$ so that 
there are two points in the inverse image of each regular point of $M$.
The so-called $2$-fold cover $M^d$ has no mirror points: 
Let $V_i$, $i=1,2, \dots$, be model neighborhoods 
covering $M$ and $V_i$ has model pairs $(U_i, G_i)$. 
The model open set $U_i$ has an induced orientation from $\bR^n$. 
For each $i$, we define a new pair of form 
$(U_i\times {\{-1, 1\}}, G_i)$ where 
$G_i$ acts by $g((x, l)) = (g(x), \sgn(g)l), l= \pm 1$ where 
$\sgn(g)$ is defined to be $1$ if $g$ is 
orientation-preserving and $-1$ if not. 
For each morphism $V_i \cap V_j \ra V_i$, we simply
define the lift
\[\tilde \phi'_{V_i\cap V_j, V_i}: V_i \cap V_j \times \{-1, 1\} 
\ra V_i \times \{-1, 1\} \]
to be $\tilde \phi_{V_i\cap V_j, V_i} \times \id_{\{-1, 1\}}$
and the homomorphism $\phi'_{V_i\cap V_j, V_i}$ to be
the old one $\phi_{V_i\cap V_j, V_i}$.
If we paste together these sets with thus-defined morphisms,
we obtain a new orbifold $M^d$. 
Also, projections $U_i \times {\{-1, 1\}} \ra U_i$
define an orbifold-covering map $M^d \ra M$ which is two-fold.
(Note that $M^d$ is not the ``doubled orbifold'' defined by Thurston
in general. For example consider the orbifold which is 
the quotient-orbifold of $\bR^3$ by the group of order-two 
generated by $-1$ times the identity map.)

\begin{lem}\label{lem:double}
Let $f:M_1 \ra M_2$ be an orbifold-map, 
and $p_1:M_1^d \ra M_1$ and $p_2: M_2^d \ra M_2$ be 
the two-fold orbifold covering maps as defined above.
Then there exists an orbifold-map $f^d: M_1^d \ra M_2^d$ 
covering $f$. 
\end{lem}
\begin{proof} 
For each model pair $(U_i \times \{ -1, 1\}, G_i)$ of $M_1^d$,
we define $f^d$ to be $f \times \id_{\{-1, 1\}}$. 
\end{proof}

\begin{prop}\label{prop:automorphisms}
Let $p:N \ra M$ be an orbifold-covering map
where $N$ and $M$ are connected orbifolds. 
Let $f:N \ra N$ be a morphism. 
Then $f$ is a covering-isomorphism.
\end{prop}
\begin{proof}
First, we consider the case when $M$ has no mirror points. 
Then we claim that $N^r$ is connected: 
In the model pair of a singular point of $N$, 
if there are no element fixing a hypersurface, 
then the set of regular points in the model open set $U$ 
is connected since the actions are conjugate to linear 
actions on sufficiently small open sets. 
Thus, each point of $U$ can be a boundary point of
only one component of $N^r$. Therefore, $N^r$ can have only
one component. (By same reason, $M^r$ is 
connected since the identity map is an orbifold-covering map.)

Now, $f|N^r: N^r \ra N^r$ is a topological covering automorphism.
Thus, $f|N^r$ is one-to-one and onto.  
By Proposition \ref{prop:nonsingular}, $f$ is a covering isomorphism.

If $M$ has mirror-points, then we form the two-fold 
orbifold-covering map $p:M^d \ra M$. Then an orbifold-map
$f^d: M^d \ra M^d$ covering $f$ is also a covering morphism 
of $p$ from $M^d$ to itself. Since $M^d$ has no mirror-points, 
$f^d$ is a covering isomorphism. Therefore, so is $f$ obviously.
\end{proof}

\begin{prop}\label{prop:modelem}
A model neighborhood of a point of an orbifold $Y$ 
is elementary for any orbifold-covering map if 
the open subset of the model pair is simply connected. 
\end{prop}
\begin{proof}
Let $V$ be a model neighborhood of $x \in Y$,
and $(\tilde V, G_V)$ the model pair, 
and $p:Y'\ra Y$ an orbifold-covering map.
For a path $f$ in $\tilde V$ with the base point 
$\tilde x$, we can lift $q\circ f$ to a path in $V'$ 
easily by using the elementary neighborhoods. 
Two homotopic paths $f$ and $f'$  
lift to homotopic path-classes again using elementary
neighborhoods. By taking a base point in $\tilde V$ and 
path-classes from the base point to all points of $\tilde V$, 
we can lift $q$ to a map $\tilde q:\tilde V \ra V'$ 
so that $p\circ \tilde q = q$. Since any path in $V'$ 
can be lifted, $\tilde q$ is obviously a surjective orbifold-map. 
(This works in the same manner as in the covering space theory.)
From here, it is straightforward to verify that 
$V'$ is of form $\tilde V$ quotient out by a finite subgroup $G'_V$
of $G_V$. That is, we show that the inverse image of 
every point of $V'$ is an orbit of $G'_V$
by an open and closedness argument.
Thus, $V$ is elementary.
\end{proof}

\begin{prop}\label{prop:ballorb}
Let $V$ be an $n$-orbifold which is a quotient space 
of an $n$-ball $\tilde V$ by a finite group $G_V$ acting on it. 
Then the following statements hold\/{\rm :}
\begin{itemize}
\item[{\rm (i)}] A connected covering orbifold $V_1$ of $V$ is isomorphic to
$\tilde V/G'_V$ for a subgroup $G'_V \subset G_V$ 
with a covering map
$p: \tilde V/G'_V \ra V=\tilde V/G$ induced  
by the identity map $\tilde V \ra \tilde V$\/{\rm ;} 
i.e., the set of the isomorphism classes of 
connected covering orbifolds
is in one-to-one correspondence with the conjugacy
classes of subgroups of $G$.
\item[{\rm (ii)}] Given two covering orbifolds 
$\tilde V/G_1$ and $\tilde V/G_2$,
a covering morphism $\tilde V/G_1 \ra \tilde V/G_2$ 
is induced by an element $g\in G:\tilde V \ra \tilde V$ 
so that $gG_1g^{-1}\subset G_2$. 
The covering morphisms are in one-to-one correspondence with
double cosets of form $G_2 g G_1$
with $g$ satisfying $gG_1 g^{-1} \subset G_2$. 
\item[{\rm (iii)}] The covering automorphism group
of a covering orbifold $V'$ is given by 
$N(G'_V)/G'_V$ where $G'_V$ is a subgroup
corresponding to $V'$ and $N(G'_V)$ is the normalizer
of $G'_V$ in $G_V$.
\end{itemize}
\end{prop}
\begin{proof}
(i) The first part follows from Proposition \ref{prop:modelem}.
The second part is a consequence of (ii).

(ii) The morphism lifts to an orbit-preserving map
$f:\tilde V \ra \tilde V$, which sends each $G_1$-orbit
to a $G_2$-orbit. Again $f$ is an element $g$ of $G_V$ 
covering the identity map of $V$. 
Thus $gG_1(x) \subset G_2(g(x))$ for 
each regular point $x$ of $\tilde V$. Thus, 
$g G_1g^{-1}\subset G_2$. 

The elements $g$ and $g'\in G_V$ induce a same map
$\tilde V/G_1 \ra \tilde V/G_2$ if and only if
$g' = g_2 g g_1$ for $g_1\in G_1$ and $g_2 \in G_2$.

(iii) follows from (ii).
\end{proof}

We first define orbifold-fiber products of orbifold-covering 
spaces of a model pairs:
Let $G_i$, $i \in I$, be a collection 
of subgroups of a finite group $G_V$
acting on an open subset $V$ of $\bR^n$. 
Then $p_i:V/G_i \ra V/G_V$ form a collection of
orbifold-covering maps.
\[p: V \times \prod_{i \in I} G_i\backslash G_V \ra V\] 
is a covering map. We let $G_V$ act on it by
\[\gamma(v, G_i \gamma_i)_{i \in I} = 
(\gamma(v), G_i \gamma_i \gamma^{-1})_{i \in I}.\]
Define the {\em orbifold-fiber-product} to be 
\[p^f: V^f 
=  (V \times \prod_{i \in I} 
G_i\backslash G_V)/G_V \ra V/G_V\] 
where $p^f$ is induced by $p$;
i.e., the orbit $[v, G_i \gamma_i]$ 
of $(v, G_i\gamma_i)$
is sent to the orbit $[v]$ of $v$.

As in Example \ref{exmp:fiberprod}, each component of 
$V^f$ is of form $V/\bigcap_{i \in I} \gamma_i G_i \gamma_i^{-1}$
for a sequence of representatives $\gamma_i$, $i \in I$, 
of cosets $G_i \backslash G_V$.  
$p^f$ is obviously given by projection, and $G_V$ 
acts transitively on the components of $V^f$.

There are covering morphisms $q_i:V^f \ra V/G_i$ 
given by $q_i:[v, G_i \gamma_i] \mapsto [\gamma_i v]$.

If we replace $\tilde V$ by $\tilde V^r$, 
we obtain an ordinary fiber product $V^{r, f}$ of
the maps $(p_i|V_i^r)_{i\in I}$. 
We can easily see that 
$V^{r, f}$ is identifiable with the inverse image
of $p^{f,-1}(V^r)$ by construction of
$V^f$. (This is nicely explained in Chapter 
5 of Thurston \cite{Thnote}. 
Also, the discussion over the regular part reduces to
Example \ref{exmp:fiberprod}.)

Let $y \in V/G_V$ be a base point and a regular point.
Given a base point $y_i$ of $V/G_i$ for each $i$ mapping 
to the base point $y$ of $V/G_V$, we can form 
a base point $y^f$ of $V^f$: Choose a base point $y_o$ in 
$V$ so that $y_i = [\gamma_i(y_o)]$ for some 
$\gamma_i \in G$. 
Then let $y = (y_o, (G_i \gamma_i)_{i\in I})$ 
and $y^f$ be the equivalence class of $y$ 
under the action of $G_V$ in $V^f$.

Again $V^f$ has a {\em universal property} that 
if $(V'', p'')$ is a connected covering space of 
$V/G_V$ with a base point $y''$ 
and $q_i:V'' \ra V/G_i$ is a covering morphism
sending $y''$ to base points $y_i$ for each $i$,
then there exists a covering morphism
$q':V'' \ra V^f$ sending $y''$ to $y^f$ 
so that $q_i\circ q' = q''_i$:
We regard $V''$ as $\tilde V/G''_V$ for 
a subgroup $G''_V$ of $G_V$. 
We can lift $q''_i$ to an orbifold-covering 
$\hat q_i:\tilde V \ra V/G_i$
for each $i$, and hence to a covering map 
$\tilde q_i:\tilde V \ra \tilde V$, which 
is given by an element 
$g_i$ of $G_V$ by Proposition \ref{prop:ballorb}
satisfying
\begin{equation}\label{eqn:gamma}
g_i G''_V g_i^{-1} \subset G_i.
\end{equation}
We map $\tilde V$ to
\[V^f = 
(\tilde V \times \prod_{i\in I} (G_i\backslash G_V))/G_V \]
by sending $u$ to the class of $(u, (G_i g_i)_{i\in I})$.
This map induces a well-defined (diagonal) morphism $q'$ 
from $V''$ to $V^J$ by equation (\ref{eqn:gamma}).
\begin{eqnarray}\label{eqn:q} 
\tilde V &\longrightarrow 
& V^f = (\tilde V \times 
\prod_{i\in I} (G_i \backslash G_V)/G_V \nonumber \\
\downarrow & & \downarrow q_i \nonumber \\
V''=\tilde V/G''_V &\stackrel{q''_i}{\longrightarrow} 
& \tilde V/G_i \nonumber \\ 
&  &\downarrow p_i \nonumber \\
&  & V = \tilde V/G_V
\end{eqnarray}
Clearly, such a diagonal map is unique.
Moreover, we may verify that $y''$ goes to $y^f$. 
Thus the component $V^f_0$ of $V^f$ containing $y^f$ has 
a universal property. Also, since we can change components, 
each component of $V^f$ or $V^f$ itself has a universal
property without base point conditions. 

Also, the universal property characterizes components of $V^f$
or $V^f$ itself up to covering isomorphisms. 
That is, if $(Z, p_Z)$ is a connected covering space of $V$ so that 
there is a covering morphism $q_{Z, i}:Z \ra V_i$,
so that $z \mapsto y_i$ for each $i$
and $Z$ satisfies the universal properties of $Y^f$ above, 
then there exists a covering isomorphism 
$L: Z \ra V^f_0$ such that $q_{Z, i} = p'_i\circ L$--(*):
The proof is obvious from the universal properties of 
$Z$ and $V^f$. 

\begin{exmp}\label{exmp:A}
Let $V$ be a closed interval and $\bZ_2$ the group of 
order two fixing a point of $V$. Let 
$p_1: V \ra V/\bZ_2$ and $p_2: V \ra V/\bZ_2$ be 
two orbifold-covering maps. 
Then the orbifold-fiber-product is obviously given 
by two copies of $V$ mapping to $V/\bZ_2$. 
\end{exmp}

Let us list a collection of orbifold-coverings
$(Y_i, p_i)$, $i\in I$, of a connected orbifold $Y$ 
for some index set $I$. Let $y_i$ be base-points of $Y_i$ 
so that $p_i(y_i) = y$ for a regular base point $y$ of $Y$.
($y_i$ are all regular.)
We now define {\em orbifold-fiber-product} of 
orbifold-covering maps $(Y_i, p_i)$:
Let us cover $Y$ by a collection $\mathcal{Z}$ of 
connected model neighborhoods  
so that the open subsets of their model pairs are simply connected.
We also require that finite intersections of model neighborhoods 
are always connected and the open subset of their model pairs 
are simply connected. Such a covering of $Y$ exists by Proposition
\ref{prop:nicecover}.

Let $V \in \mathcal{Z}$ 
be a model-neighborhood of $Y$ with a model pair
$(\tilde V, G_V)$ where $\tilde V$ is simply-connected. 
Take a component $V_i^j$, 
of $p_i^{-1}(V)$ for each $Y_i$. 
Let $G_i^j$ denote the subgroup of $G_V$ 
so that $\tilde V \ra V_i^j$ is the covering 
map given as a quotient map for $G_i^j$ 
acting on $\tilde V$.  
For each choice of $j$ for $i$, we define 
a map $J$ defined on $I$ 
to the set of components of $p^{-1}_i(V)$ for $i \in I$
so that $J(i)$, also denoted by $j(i)$,
is a component of $p^{-1}_i(V)$.
With a fixed function $J$,
we form a fiber product $V^J$ of 
$V_i^{j(i)}, i\in I$.
We define
\[V^J = (\tilde V \times 
\prod_{i\in I} (G_i^{j(i)}\backslash G_V))/G_V\]
where the $G_V$-action is given by 
\[\gamma(v, (G_{i}^{j(i)}\gamma_{i})_{i\in I})
= (\gamma v, (G_i^{j(i)}\gamma_{i}\gamma^{-1})_{i\in I}).\]

We define $q_J:V^J \ra V$ to be the obvious covering map
sending $(u, (G_i^{j(i)}\gamma_i)_{i\in I})$ to the class of $u$ 
in $V$, i.e., the orbit $[u]$ of $u$. 
From $V^J$, we can define a covering map
$q_i: V^J \ra V_i^{j(i)}$ by 
sending $[u, (G_i^{j(i)} \gamma_i)_{i\in I}]$
to $[\gamma_i u]$ in the equivalence class of
$\tilde V/G_i^{j(i)}$. This clearly is
a well-defined orbifold-morphism.

We note the universal property of $V^J$ 
that given a sequence of morphisms
$q''_i:V'' \ra V_i^{j(i)}$ for all $i$, 
there is a covering morphism
$q':V'' \ra V^J$ so that $q_i\circ q' = q''_i$:

Now, we define $\hat V$ as the disjoint union 
$\coprod_J V^J$ for all functions $J$.
It has an obvious covering map
$\hat p:\hat V \ra V$. We can define a morphism 
\[q_i:\hat V \ra \bigcup_{j} V_i^j = p_i^{-1}(V)\] 
by defining $q_i$ as above for each of $V^J$. 

Since $\coprod_J V^J$ contains the fiber products of 
\[p_i| p_i^{-1}(V^r): p_i^{-1}(V^r) \ra V^r,\]
we see that if the base point $y$ of $Y$ is in $V$, then there exists 
a regular point $y^f$ mapping to $y_i$ under $q_i$ for 
each $i$. Thus, we construct $y^f$ in this manner, later 
to be identified.

Also, $\hat V$ has a following universal property: 
given a sequence of morphisms 
\[q''_i:V'' \ra \coprod_j V_i^j = p_i^{-1}(V)\] 
for all $i$ so that $p_i\circ q''_i$ is a fixed 
covering map $V'' \ra V$, there exists 
a unique morphism $q':V'' \ra \hat V$ so
that $q_i\circ q' = q''_i$. This follows from considering 
where each component of $V''$ maps to. 
The covering $\hat V$ is said to be a {\em fiber product}
of $p_i^{-1}(V), i\in I$.

Let $U$ be a connected open subset of $V$, such 
as $U = V \cap V'$ for another model neighborhood $V'$
in the covering $\mathcal{Z}$.
We assume that $U$ is modeled on a pair
$(\tilde U, G_U)$. Then components of
$q^{-1}(U)$ in $\tilde V$ are homeomorphic 
to $\tilde U$ and the subgroup of
$G_V$ acting on a component is 
isomorphic to $G_U$. 

Let $U_i^{j, k}$ denote the components of
$p^{-1}(U)$ in $V_i^j$ for each $i, j$. 
Let $G_{U, i}^{j, k}$ denote a subgroup of
$G_U$ so that $\tilde U/G_{U, i}^{j, k}$
equals the covering $U_i^{j, k} \ra U$.

Let $K$ be a function defined on $I$ 
by sending $i$ to an index $k(i)$ among the indices $k$ of
components of form $U_i^{j(i), k}$.
Let the fiber product 
\[U^{J, K} = (\tilde U \times 
\prod_{i\in I} (G_{U, i}^{j(i), k(i)}\backslash G_U))
/G_U\]
be defined where $G_U$ acts by 
\[(u, (G_{U, i}^{j(i), k(i)}\gamma_i)_{i\in I}) \mapsto
(\gamma u, (G_{U, i}^{j(i), k(i)}\gamma_i\gamma^{-1})_{i\in I}),
\hbox{ for } \gamma \in G_U.\] 
Let 
\[q_{U, i}^{J, K}:U^{J,K} \ra U_{U, i}^{j(i), k(i)} \subset p_i^{-1}(U)\]
be the morphism defined by sending 
$(u, (G_i^{j(i), k(i)}\gamma_i)_{i\in I})$ to 
the equivalence class of $\gamma_i u$.
Let $p_U^{J, K}:U^{J, K} \ra U$ denote the covering map.

Define $U^J$ by taking the disjoint union 
$\coprod_{K} U^{J, K}$, and $\hat U$
by $\coprod_J U^J$. We define 
morphisms 
\[q_{U, i}^J: U^J \ra \coprod_k U_i^{j(i), k} \subset p_i^{-1}(U)\]
by restricting it to be $q_{U, i}^{J, K}$ for appropriate 
components, and define 
morphisms 
\[q_{U, i}: \hat U \ra \coprod_{j, k} U_i^{j, k} = p_i^{-1}(U)\] 
similarly.
We let $p_U^J: U^J \ra U$ and
$\hat p_U: \hat U \ra U$ denote the covering maps. 
We note that $\hat U$ has the appropriate universal property also:
i.e., if 
\[q''_i: U'' \ra \coprod_{j, k} U_i^{j, k}= p_i^{-1}(U)\] 
is a morphism for each $i$, then there exists a unique morphism
$q''_U:U'' \ra \hat U$, so that 
\begin{equation} \label{eqn:**}
q_{U, i} \circ q''_U = q''_i. 
\end{equation}

We will now identify $\hat p^{-1}(U)$ in $\hat V$ with $\hat U$:
Since $q_i: \hat V \ra p_i^{-1}(V)$ is a morphism, 
$\hat p^{-1}(U)$ is mapped to $p_i^{-1}(U)$ by $q_i$. 
Thus, there is a morphism $f: \hat p^{-1}(U) \ra \hat U$
by the universal property of $\hat U$.
By construction, 
$f$ sends $(u, (G^{j(i)}_i \gamma_i)_{i\in I})$
for $p(u) \in U$ to a point of $\hat U$ mapping 
to $q_i(\gamma_i(u))$ under $q_{U, i}$ for each $i$.  
We obtain a morphism
\[f|\hat p^{-1}(U^r):\hat p^{-1}(U^r) \ra \hat p_U^{-1}(U^r),\]
which is an ordinary covering-isomorphism between
fiber products of ordinary covering spaces  
$\hat p^{-1}(U^r)$ and $\hat p_U^{-1}(U^r)$ of $U^r$.
(This follows since the two sets are obviously topological
fiber-products over $U^r$ and $f$ is the natural identification
with a commutative diagram: 
\begin{eqnarray*}
\hat p^{-1}(U^r)  & \stackrel{f}{\longrightarrow} & \hat p^{-1}_U(U^r) \\
q_i \downarrow & & q_{U, i} \\
p_i^{-1}(U^r) & \stackrel{\id}{\longrightarrow} & p_i^{-1}(U^r).)
\end{eqnarray*}
By Proposition \ref{prop:nonsingular}, 
we can identify $\hat U$ as a suborbifold of $\hat V$. 
Since $q_{U, i} \circ f = q_i$ by the uniqueness part of 
equation \ref{eqn:**}, the orbifold-map
$q_i$ on $\hat V$ extends $q_{U, i}$ on the suborbifold $\hat U$. 
Also, since $\hat p_U \circ f = p$ while $f$ is a morphism, 
$\hat p_U$ is extended to $\hat p:\hat V \ra V$. 

We let $\hat Y^f$ be the quotient space of 
union of all $\hat U$s as $U$ ranges over all open model 
subsets in the covering $\mathcal{Z}$ 
with identification given as above. 
We call $\hat Y^f$ the {\em orbifold-fiber-product} of $Y_i$s. 
We obviously have an orbifold-map $p^f: \hat Y^f \ra Y$ extended from 
$\hat p$s defined over the model neighborhoods.
Obviously, base points $y^f$s in $\hat V$s 
correspond to a unique point in $Y^f$, 
to be denoted by $y^f$ again.

We choose a set $\hat Y$ as a component of $\hat Y^f$
containing the base point $y^f$. Let $\hat p: \hat Y \ra Y$ 
be the restriction of $p^f$.
 
The topology of $\hat Y$ is given by the basis which 
are sets of form of components of $V^J$ as 
$V$ ranges over the elementary neighborhoods of $Y$.
There are well defined morphisms $\hat p: \hat Y \ra Y$ 
and $q_i: \hat Y \ra Y_i$ extending $p$ and $q_i$ on
each sets of form $\hat V$. The extension exists 
by equation \ref{eqn:**}. We say that $\hat Y$ with $\hat p$ is 
a {\em component of a fiber product} of $Y_i$s. 

We show that $\hat Y$ is Hausdorff and second-countable:
Let $x$ and $y$ be two points of $\hat Y$. 
If $\hat p(x)$ and $\hat p(y)$ are distinct, then
choose disjoint elementary neighborhoods in $Y$ 
and their inverse images are disjoint open sets. 
If $\hat p(x)$ equals $\hat p(y)$, then $x$ and 
$y$ are in different sets of form of a component of $V^J$ for some 
elementary open set $V$, then as components of $V^J$ are basis
elements, $x$ and $y$ are contained in disjoint 
open sets. If $x$ and $y$ are in the same component $U$ of 
some $V^J$, then the inverse image in $\tilde V$
of $x$ and that of $y$ meet up to an action
of a finite group $G_V$ where $(\tilde V, G_V)$ 
is the model pair for $V$ considering $U$ 
as a quotient space of $\tilde V$. Since $x$ and $y$ are 
distinct, there are certainly disjoint neighborhoods 
in $U$ as $U$ is a quotient space of $\tilde V$ by 
a proper subgroup of $G_V$ where $x$ and $y$
are in different orbits. 

Since each open set of form $V^J$ is locally compact,
$\hat Y$ is locally compact, and $\hat Y$ is metrizable.
Since a component of the topological fiber product $\hat Y^r$ 
of $Y_i^r$ is separable,
so is $\hat Y$; thus, $\hat Y$ is a Hausdorff second countable set. 
Since $\hat Y$ is covered by model neighborhoods,
$\hat Y$ is an orbifold, and $\hat p:\hat Y \ra Y$ is
an orbifold-covering map. 

Similarly, components of $\hat Y^f$ enjoy 
the same properties, i.e., they are Hausdorff second countable 
and hence are orbifolds. Thus, 
the fiber-product $\hat Y^f$ is a disjoint union of orbifolds.
Therefore, $p^f: \hat Y^f \ra Y$ is an orbifold-covering map,
i.e., the aim of our construction.

We verify a universal property: 
Let $(Z, p_Z)$ be a covering orbifold of $X$
so that there are covering morphisms 
$q^Z_i:Z \ra Y_i$ sending the base point $z$ to $y_i$s
where $p_Z(z) = p_i(y_i) = y$.
We show that there exists 
a covering morphism $Z \ra \hat Y^f$ sending $z$ to $y^f$:
Define $Z^r$ to be the inverse image of 
$Y^r$ in $Y$. Then $Z^r$ is an open dense set,
and it has a morphism to $\hat Y^{f, r}$ over $X^r$
sending $z$ to $y^f$
since $\hat Y^{f,r}$ is a fiber-product of $Y_i^r$s.
Now, there is a unique extension of this map
from $Z$ to $\hat Y^f$ since the extensions are given 
by the universal properties of fiber-products of 
the inverse images of model neighborhoods.

\begin{exmp}\label{exmp:thurston}
For reader's convenience, let us explain 
Thurston's example in Chapter 5 of \cite{Thnote}. 
We consider two orbifold-coverings of an interval 
$I=[-1,1]$ with two mirror points at the end. 
Then $I$ is covered by 
a circle $\SI^1$ where the covering map is given by
a projection to the $x$-axis. This is a regular covering 
with $\bZ_2$ acting by reflections in 
the $x$-axis. Let $p_1: \SI^1 \ra I$ be this covering. 
The next covering is an orbifold-covering $p_2: J \ra I$ where 
$J$ is an interval $[-1, 3]$ with two mirror points at the end.
Define $p_2: J \ra I$ by $p_2(t) = t$ if 
$-1 \leq t \leq 1$ and $p(t) = -t + 2$ if $1 \leq t \leq 3$. 
Again, this is a regular covering.

Cover $I$ by three open sets $I_1 = [-1, -\eps), 
I_2 = (-2\eps, 2\eps), I_3 = (\eps, 1]$ where 
$0 < \eps < 1/2$. The above construction 
tells us that the fiber product of inverse images of 
$I_2$ is a union of four arcs. 
Over $I_1$, it is a union of two arcs. 
Over $I_3$, it is a union of
two arcs. Since the fiber products over $(-2\eps, -\eps)$ and 
$(\eps, 2\eps)$ are unions of four arcs obviously 
identifiable as subsets, we clearly see that the fiber product 
is isomorphic to a four-fold covering $\SI^1 \ra I$ 
with four open arcs mapping to $I_1$ or $I_3$ as two-fold 
coverings. 
\end{exmp}

In the note of Thurston \cite{Thnote}, he proved 
that each orbifold $Y$ has a so-called universal covering orbifold
$\tilde Y$ with an orbifold-map $p_Y:\tilde Y \ra Y$
so that given an orbifold-covering map $p_Z:Z \ra Y$
where $Z$ has a connected underlying space, 
there is an orbifold-map $q:\tilde Y \ra Z$ so that 
$p\circ q$ equals $p_Y$. 

We give a more precise definition:
Let $y$ be a regular base point of $Y$.
A {\em universal covering orbifold} is 
a connected covering $(\tilde Y, p_Y)$ with a regular base point 
$\tilde y$ mapping to $y$ so that for any covering $p_Z:Z \ra Y$
where $Z$ is a connected orbifold with a regular base point $z$ over $y$
there is an orbifold-covering-morphism 
$q: \tilde Y \ra Z$ so that $q(\tilde y) = z$ and 
$p_Z \circ q = p_Y$. We require that this should holds for any choices 
of regular base points $y, \tilde y \in p_Y^{-1}(y)$, 
and $z \in p_Z^{-1}(y)$.

The uniqueness of the universal covering of an orbifold 
up to covering isomorphisms (preserving base points) is obvious
from the definition.

The group of automorphisms of a universal cover of
$Y$ is said to be
the {\em fundamental group}, and we denote it by $\pi_1(Y)$. 

\begin{prop}[Thurston]\label{prop:univ}
Let $Y$ be a connected orbifold. Then
there exists a universal covering orbifold $\tilde Y$
unique up to covering isomorphism.
Moreover, the fundamental group of $\tilde Y$ acts 
transitively on the inverse image of the base point $y$.
\end{prop} 
\begin{proof}
Let $\hat Y$ be obtained by a list of coverings 
$(Y_i, p_i)$ with base point $y_i$ over a base point $y$
which lists one-element from each
isomorphism class of covering maps of $Y$ 
preserving base-points. 
Let $y^f$ be the base point from the above construction.

Given any covering $p:Z \ra Y$ with a base point $z$, 
since we obviously have our $Z$ isomorphic to say $Y_i$, 
$i\in I$, there exists a covering morphism
$q_Y:\hat Y \ra Z$ where $q_Y(y^f) = z$. 

Let $y'$ be a point of $p^{f, -1}(y)$ different from $y$. 
We show that there exists a deck transformation sending 
$y^f$ to $y'$: Clearly, $(\hat Y, \hat p)$ with $y'$ as 
a base point is in the list of all covering maps of $Y$. 
Thus, there exists a morphism $g: \hat Y \ra \hat Y$ 
sending $\tilde y$ to $y'$. By Proposition \ref{prop:automorphisms}, 
$g$ is a deck transformation.

Now, let $(Z, p_Z)$ be a cover of $Y$ with a base point $z$
mapping to a base point $x$ of $Y$, perhaps different from above $y$. 
Let $y'$ be a point of $p^{f, -1}(x)$ of $\hat Y$ (deemed to be 
our new base point). 
Find a path $\alpha$ from $y'$ to $\tilde y$, which 
maps to a path $\alpha'$ from $x$ to $y$ on $Y$, 
and find a path $\alpha''$ on $Z$ from $z$ ending at a point $z'$
lifting $\alpha'$. 
Then $p_Z(z') = y = p^f(\tilde y)$. 
By above construction of fiber-products, 
there exists a morphism $g:\hat Y \ra Z$ 
so that $g(y'') = z'$ for some $y'' \in \hat p^{-1}(y)$
since $\hat Y$ and $Z$ are connected.
By precomposing with $g$ a deck transformation, we may assume
that $g(\tilde y) = z'$. 
Since $\alpha$ goes to $\alpha''$, 
we see that $g(y')$ maps to $z$. 
Therefore $\hat Y$ satisfies the definition of a universal cover.
\end{proof}

\begin{prop}\label{prop:univlift} 
Let $p_1:\tilde Y \ra Y_1$ and $p_2: \tilde Y \ra Y_2$ 
be universal covering orbifold-maps. 
Then any orbifold-map $f: Y_1 \ra Y_2$ 
covering an orbifold-diffeomorphism $g: Y \ra Y$ lifts 
to an orbifold-diffeomorphism $\tilde f:\tilde Y \ra \tilde Y$.
The lift is unique if we decide the 
value $\tilde f(y_0)$ among $p_2^{-1}(f(p_1(y_0)))$
for a base-point $y_0$ of $\tilde Y$, where
$\tilde f(y_0)$ can be chosen to be any such point.
Finally, if $g$ is the identity, then $\tilde f$ is 
an automorphism.
\end{prop}
\begin{proof}
Since $f: Y_1 \ra Y_2$ is a covering map, 
$f\circ p_1$ is the universal covering of $Y_2$. 
Since $p_2$ is also a universal covering of $Y_2$, 
there is a morphism $\tilde f: \tilde Y \ra \tilde Y$ 
so that $p_2 \circ \tilde f = f \circ p_1$ by the uniqueness
of the isomorphism class of universal covering orbifolds.  
The uniqueness follows since $\tilde f$ restricts to
$\tilde Y^r$ a lift of $f$ restricted $Y_1^r$
and the ordinary covering space theory.
The freedom of choice follows by the transitivity of
automorphism group in the fiber of $\tilde Y \ra Y_2$.
The final statement follows trivially.
\end{proof}

The fundamental group $\pi_1(Y)$ acts on 
the universal cover $\tilde Y$
as a group of orbifold-diffeomorphisms. The action is proper 
since given a point of a model neighborhood of $\tilde Y$
which is a component of a preimage of a model neighborhood of $Y$, 
it may return to the neighborhood only finitely many times 
since the only equivalent points are in the orbits. 
Therefore, $\tilde Y/\pi_1(Y)$ is again an orbifold. 

We obtain most of the useful results of the covering space 
theory for topological spaces in the orbifold case.
\begin{cor}\label{cor:univ}
\begin{itemize}
\item[{\rm (i)}] The fundamental group acts transitively
on each fiber of a universal cover $\tilde Y$, 
i.e, on the inverse image of a regular point $x$
of $Y$. Moreover, $\tilde Y/\pi_1(Y) = Y$.
\item[{\rm (ii)}] Each covering space $p_1:Y_1 \ra Y$ is isomorphic 
to a covering map $p': \tilde Y/\Gamma \ra Y$ 
where $p'$ is induced from the universal covering map 
$p:\tilde Y \ra Y$ and $\Gamma$ is a subgroup of $\pi_1(Y)$.
\item[{\rm (iii)}] The isomorphism classes of covering spaces of $Y$ 
are in one-to-one correspondence with the conjugacy classes 
of subgroups of $\pi_1(Y)$.
\item[{\rm (iv)}] The group of automorphisms of a covering 
space $\tilde Y/\Gamma$ is isomorphic to $N(\Gamma)/\Gamma$
where $N(\Gamma)$ is the normalizer of $\Gamma$ in $\pi_1(Y)$.
\item[{\rm (v)}] A covering $Y' \ra Y$ is regular if and only if 
$Y'$ is isomorphic to $\tilde Y/\Gamma$ for a normal
subgroup $\Gamma$ of $\pi_1(Y)$.
\item[{\rm (vi)}] Let $Y_1$ and $Y_2$ be orbifolds with 
universal covering orbifolds $\tilde Y_1$ and $\tilde Y_2$ 
respectively. A lift $f: \tilde Y_1 \ra \tilde Y_2$ of 
an orbifold-diffeomorphism $g: Y_1 \ra Y_2$ is an isomorphism.
\end{itemize}
\end{cor}
\begin{proof}
(i) This follows from definition, i.e., the condition
on base-points. The covering map $p_Y: \tilde Y \ra Y$ 
induces an orbifold-map 
\[\tilde Y/\pi_1(Y) \ra Y\]
which is one-to-one over the regular part. 
$p_Y$ restricts to a homeomorphism
over the regular part since it is proper 
and is a local homeomorphism.
Therefore, $p_Y$ is an orbifold-diffeomorphism
by Proposition \ref{prop:nonsingular}.

(ii) There is a morphism $q:\tilde Y \ra Y_1$. 
Since $\tilde Y$ covers any cover of $Y_1$ also, 
$\tilde Y$ is again a universal cover of $Y_1$. 
Since $Y_1$ is a quotient of its universal cover, 
(ii) follows.

(iii) If two covering spaces $p_1:\tilde Y/\Gamma_1 \ra Y$ and 
$p_2:\tilde Y/\Gamma_2 \ra Y$ are isomorphic, then there exists 
a morphism $f: \tilde Y/\Gamma_1 \ra \tilde Y/\Gamma_2$. 
$f$ lifts to a morphism $\gamma: \tilde Y \ra \tilde Y$ by Proposition 
\ref{prop:univlift}. Since $f$ is a morphism,
$f$ satisfies $p_2 \circ f = p_1$, and so $f$ covers identity on $Y$.
Thus, the lift $\gamma$ of $f$ is a deck transformation 
of $p_Y$. In order that $\gamma$ descends to a map $f$, 
we need that for each $\alpha \in \Gamma_1$, 
there exists $\alpha' \in \Gamma_2$ so that 
$\gamma \alpha = \alpha' \gamma$. Thus, 
$\gamma \Gamma_1 \gamma^{-1} \subset \Gamma_2$. 
A converse argument shows that 
a conjugate of $\Gamma_2$ is in $\Gamma_1$. 
Hence, $\Gamma_1$ and $\Gamma_2$ are conjugate. 

(iv) This follows from (iii).

(v) Given a base point $y'$ of $\tilde Y$ mapping 
to a base point $y$ of $Y$, each point of $p_Y^{-1}(y)$ is 
of form $\gamma(y')$ for $\gamma \in \pi_1(Y)$. 
Given $\tilde Y/\Gamma$ for a subgroup $\Gamma$ of $\pi_1(Y)$, 
and two points in the fiber of $y$, we see that
there exists a deck transformation $\gamma$ 
sending $[y']$ to $[\gamma'(y')]$ for some $\gamma'$
if and only if $\gamma$ normalizes $\Gamma$. 
It is easy to see that all coset representative of 
$\pi_1(Y)/\Gamma$ has to occur. Thus $\Gamma$ is normal. 

(vi) Let $p_1:\tilde Y_1 \ra Y_1$ and 
$p_2:\tilde Y_2 \ra Y_2$ be orbifold-covering maps. 
Then we have $p_2 \circ f = g \circ p_1$. The map
$p_2\circ f$ is a universal covering map of $Y_2$ since 
$\tilde Y_1$ can cover any orbifold-covering $Y_2$ using $f$
and deck transformation group of $Y_2$. 
Thus, there exists an isomorphism $h: \tilde Y_1 \ra \tilde Y_2$.
If we identify $Y_1$ and $Y_2$ by $g$ and $\tilde Y_1$ and 
$\tilde Y_2$ by $h$, then $f$ becomes a morphism. 
By Proposition \ref{prop:automorphisms}, 
$f$ is an isomorphism.
\end{proof}

\begin{prop}\label{prop:lifthomt}
Let $Y$ be an orbifold and $f: Y \ra Y$ and $g:Y \ra Y$ 
orbifold-diffeomorphisms with homotopy $H: Y \times I \ra Y$. 
Then for any choice of lift $\tilde f: \tilde Y \ra \tilde Y$ 
of $f$, there is a unique lift $\tilde H: \tilde Y \times I \ra \tilde Y$
which becomes a homotopy between $\tilde f$ and a lift 
$\tilde g: \tilde Y \ra \tilde Y$ of $g$.
\end{prop}
\begin{proof}
Clearly, the lift $H': \tilde Y^r \times I \ra \tilde Y^r$ of 
$H| Y^r \times I: Y^r \times I \ra Y^r$, which is a homotopy between 
$\tilde f| Y^r: Y^r \ra Y^r$ and a map $g': \tilde Y^r \ra \tilde Y^r$, 
exists. We give a local description. 

Let $(x, t)$ be a point of $Y\times I$. There is a model open subset
$(N\times J, G_N)$ of $(x, t)$ where $G_N$ acts so 
that $(y, s) \mapsto (g(y), s)$ for $g \in G_N$
where $J$ is a small open interval containing $t$.
$H$ lifts to $H^l: N \times J \ra M$ where $M$ is 
a model open subset where a finite group $G_M$ also acts on, 
and there is a homomorphism $k:G_N \ra G_M$
so that $H^l(g(y), s) = k(g)(H^l(y, s))$.
We also have $f(g(y)) = k(g)(f(y))$ for $g \in G_N$ and 
$k(g) \in G_M$.

Let $q(N\times J)$ be the neighborhood of $(x, t)$ in $Y \times I$
corresponding to $N\times J$.
Take a component $N'\times J$ of the inverse image of $q(N\times J)$ in 
$\tilde Y \times I$. Then $\tilde f$ maps it into a component $M'$ of 
the inverse image of $M$ in $\tilde Y$. Moreover, we may 
assume without loss of generality 
that $\tilde f$ has the same model as above $N$ and $M$ 
and an associated homomorphism $k: G'' \ra G'''$
where $G''$ is a subgroup of $G_N$ and 
$G'''$ one of $G_M$.

The map $H'$ is defined on $N^{\prime, r} \times J$ 
mapping to $M^{\prime, r}$ and it extends $\tilde f|N^{\prime, r}$.
Clearly, $H'|N^{\prime, r} \times I$ lifts to $H^l$ on 
$N^r\times J \ra M^r$ with an associated homomorphism $k$ restricted 
to $G'' \ra G'''$. This means that $H^l$ induces 
an orbifold-map $H''| N'\times J: N'\times J \ra M'$. 
Since $H''$ extends $H'$ in the neighborhood,
by considering every neighborhood, 
we can show that $H'$ extends to a homotopy 
$\tilde H: \tilde Y \times I \ra \tilde Y$. Obviously,
a lift $\tilde g$ is obtained by restricting $\tilde H$.
\end{proof}

Given orbifolds $M$ and $N$, and an 
orbifold-diffeomorphism $f:M \ra N$
which lifts to a diffeomorphism $\tilde f:\tilde M \ra \tilde N$,
we obtain
an induced homomorphism $\tilde f_*: \pi_1(M) \ra \pi_1(N)$:
For each deck-transformation
$\vth$ of $\tilde M$, let $\tilde f_*(\vth)$ be 
the deck-transformation 
$\tilde f \circ \vth \circ \tilde f^{-1}$.

By Proposition \ref{prop:lifthomt}, given orbifold-diffeomorphisms 
$f_1, f_2: M \ra N$ with a homotopy $H$, a lift
$\tilde f_1:\tilde M \ra \tilde N$ of $f_1$ is homotopic to
a lift $\tilde f_2:\tilde M \ra \tilde N$ of $f_2$ by a lift 
of a homotopy $\tilde H: \tilde M \times [0,1] \ra \tilde N$. 

\begin{prop}\label{prop:diffeomorphisms}
If $\tilde f_2: \tilde M \ra \tilde N$ is 
a diffeomorphism homotopic to $\tilde f_1$ 
by a homotopy $h:\tilde M \times [0,1] \ra \tilde N$ 
equivariant with respect to $\tilde f_{1*}:\pi_1(M) \ra \pi_1(N)$,
then $\tilde f_{2 *} = \tilde f_{1 *}$.
\end{prop}
\begin{proof}
Let $\gamma$ be a deck transformation of $\tilde M$.
We have that $\gamma'=\tilde f_2\circ \gamma \circ \tilde f_2^{-1}$
is homotopic to $\gamma''=\tilde f_1\circ \gamma \circ \tilde f_1^{-1}$,
and let $H$ be the homotopy between them
given by $H_t = h_t \circ \gamma \circ h_t^{-1}$ for $t\in [0,1]$.
Then $H_t:\tilde N \ra \tilde N$ is a deck transformation
for each $t \in [0,1]$. (To see this simply post-compose $H_t$ with the 
covering map of $N$.)
Since the group of deck transformations
is discrete in $C^s$-topology,
$\gamma'$ and $\gamma''$ are equal. 
\end{proof}

\begin{rem} 
We are not aware of the full theory of
liftings of maps of orbifold-covering spaces using fundamental groups.
But it might be desirable to have one for other purposes than 
required in this paper. At any rate, the inclusion of
such a theory will fully complete the orbifold-covering-space theory,
which might be an interesting project. 
\end{rem}

\begin{rem} For two-dimensional orbifolds,
the constructions of the universal covers are considerably
simpler, and are exposed in Scott \cite{Scott:83}.
\end{rem}

\begin{rem}\label{rem:good}
A {\em good} orbifold is an orbifold with 
a covering orbifold that is a manifold. 
Since its universal covering orbifold covers a manifold, 
each of the model pairs of the universal covering orbifold 
has a trivial group action. Thus, the universal covering 
orbifold is a manifold. 

A {\em very good} orbifold is an orbifold with 
a finite regular cover that is a manifold.
 
A good orbifold $Y$ is always orbifold-diffeomorphic
to $M/\Gamma$ where $M$ is a simply-connected 
manifold and $\Gamma$
is a discrete group acting on $M$ properly.

A good orbifold $M$ has a covering that is 
a simply-connected manifold $\tilde M$. Then it is a universal
covering orbifold.
Finally, if $Y = M/\Gamma$ and $M$ is simply connected, 
then $\pi_1(Y)$ equals $\Gamma$. 
\end{rem}

\section{$(G,X)$-structures on orbifolds}

A {\em $(G, X)$-structure} on an orbifold $M$ is a collection of
charts $\phi_U:U \ra X$ for each model pair $(U, H_U)$ 
so that $\phi_U$ conjugates the action of $G_U$ with
that of a finite subgroup $G_U$ of $G$ on $\phi(U)$ 
by an isomorphism $i_U:H_U \ra G_U$, 
and the inclusion map induced 
map $U \ra V$ is always realized by an element $\vth$ of $G$
and the homomorphism $G_U \ra G_V$ is given by a conjugation by $\vth$;
i.e., $g \mapsto \vth \circ g \circ \vth^{-1}$.
A different choice equals $\vpi\vth$ for $\vpi \in G_U$ 
and the homomorphism change by conjugation by $\vpi$.  
(Again, we need the assumption that the action of $G$ is 
locally faithful)

A maximal such family of collections $(\phi_U, i_U)$ is said 
to be a {\em $(G, X)$-structure} of $M$. 
A $(G,X)$-structure on $M$ induces $(G,X)$-structures on
its covering orbifolds.

If an orbifold $M$ has a $(G, X)$-structure, then 
we can choose a model pair $(U, G_U)$ for each model
set where $U \subset X$ and a subgroup $G_U$ of $G$
using the charts. 

A {\em $(G,X)$-map} $f$ between two $(G,X)$-orbifolds $M$ and $N$ is 
a map so that for each point $x$ of $M$ and a point $y$ of $N$ 
so that $x = f(y)$, and a neighborhood $U$ of $x$ modeled on
a pair $(\tilde U, H_U)$ with a chart $\phi_U$ and 
an isomorphism $i_U: H_U \ra G_U \subset G$, 
there is a neighborhood $V$ of $y$ modeled on
a pair $(\tilde V, H_V)$ with a chart $\phi_V$ and 
an isomorphism $i_V: H_V \ra G_V$ so that 
$f$ lifts to a map $\tilde f:\tilde V \ra \tilde U$ 
equivariant with respect to a homomorphism $H_V \ra H_U$ 
induced by a homomorphism $G_V \ra G_U$ given by
a conjugation $g \mapsto \vth g \vth^{-1}$ by some $\vth \in G$. 

\begin{thm}[Thurston] 
A $(G, X)$-orbifold $M$ is a good orbifold.
There exists an immersion $D$ from the universal 
covering manifold 
$\tilde M$ to $X$ so that 
\[D\circ \vth = h(\vth)\circ D, \vth\in \pi_1(M)\]
holds for a homomorphism $h:\pi_1(M) \ra G$,
where $D$ is a local $(G,X)$-map. Moreover, any such 
immersion equals $g\circ D$ for $g\in G$, 
with the associated homomorphism $g\circ h(\cdot)\circ g^{-1}$.
\end{thm}
\begin{proof}
This is found in Chapter 5 of Thurston \cite{Thnote}
(see also Bridson-Haefliger \cite{Hae:99} and 
Matsumoto and Montesinos-Amilibia \cite{MatMont}).
We rewrite it here for the reader's convenience: 
Let $N$ be a neighborhood of $x \in \Sigma$,
and $(\tilde N, H_N)$ the model pair
for $\tilde N$ an open set in $X$ and $H_N$ 
the associated finite group acting on $\tilde N$.
We form 
$G\times \tilde N$ and give an action of $H_N$ 
by $\gamma(g, y) = (\gamma g, \gamma y)$. 
Then $G(N)= (G \times \tilde N)/H_N$ is a manifold 
and has a projection $p_N:G(N) \ra N$ induced 
by the projection to the second factor. 
(Here, $p_N$ is an orbifold-map.)

We find a nice locally finite cover of $M$ by 
such neighborhoods $\{N_1, N_2, \dots\}$.
If $N_i$ and $N_j$ meet, then $N_i \cap N_j$ 
has an inclusion map $i: N_i\cap N_j \ra N_i$. 
Then there is a connected open subset $A$ of $\tilde N_i$ 
and a subgroup $H_A$ acting on it being 
a model for $N_i\cap N_j$.
We form $G(N_i\cap N_j)_A$ where $A$ denotes
the fact we used $A$ as a model and find a map
$\tilde i: G(N_i\cap N_j)_A \ra G(N_i)$ 
induced by $G\times A \ra G\times \tilde N_i$ defined by
$(g, x) \mapsto (g, x)$. 
The map $\tilde i$ is a well-defined imbedding since 
a different choice $i$ gives us $\vth \circ i$, $\vth \in G$,
and so $\tilde i$ is replaced by a map 
$(g, x) \mapsto (\vth g, \vth x)$.

We find an open subset
$B$ of $\tilde N_j$ corresponding to $N_i\cap N_j$,
and form $G(N_i \cap N_j)_B$ similarly,
and find an imbedding $G(N_i\cap N_j)_B \ra G(N_j)$.
Since $G(N_i\cap N_j)_A$ and $G(N_i\cap N_j)_B$
can be identified by the identification of the 
model pairs, we see that $G(N_i)$ and $G(N_j)$ 
can be pasted by this relation. 
We can easily show that such identifications 
of $G(N_1), G(N_2),\dots $ are possible,
and obtain a manifold $G(M)$ from the identifications.

The foliation of $G(N_i)$ with leaves 
that are images of $g \times \tilde N_i$ for $g\in G$
gives rise to a foliation on $G(M)$ whose leaves
meet the fibers of $p_N$ at unique points. 
Take a leaf $L$ in $G(M)$. Then $p_N|L: L \ra M$
is an orbifold-covering map and $L$ is a manifold. 
Take a universal cover $\tilde L$ of $L$ with covering 
map $p_L$. Then $p_N \circ p_L$ is a universal covering 
map of $M$. $L$ has a $(G,X)$-structure since 
it covers $M$: one can induce charts. 
Then $\tilde L$ has a $(G,X)$-structure.

Since by Remark \ref{rem:good} 
$\tilde L$ is a universal cover of $M$, Corollary \ref{cor:univ}
implies that $\tilde L/\Gamma$ for the fundamental group $\Gamma$ is 
$(G,X)$-diffeomorphic to $M$ by a map induced by $p_N\circ p_L$.
As $\tilde L$ is a $(G,X)$-manifold, $M$ has a developing map
$D:\tilde L \ra X$ (which follows from the geometric structure 
theory for manifolds).
For a deck transformation $\gamma$,
$D\circ \gamma$ is also a $(G,X)$-map, and this means
that $D\circ \gamma = h(\gamma) \circ D$ for 
some $h(\gamma) \in G$. We can clearly verify that 
$h:\Gamma \ra G$ is a homomorphism.
The rest of the conclusion follows in the same way 
as the geometric structure theory for manifolds. 
\end{proof}

\begin{rem}
In most cases, geometric orbifolds are also very good
due to Selberg's lemma since our Lie groups are often subgroups of 
linear groups.
\end{rem}

We now assume that $M$ is a compact $(G, X)$-orbifold with 
a universal cover $\tilde M$.
A pair $(D, h)$ of immersions $D:\tilde M \ra X$ equivariant
with respect to a homomorphism $h:\pi_1(M) \ra G$ is 
said to be a development pair of $M$. $D$ is called
a {\em developing map} and $h$ a {\em holonomy homomorphism}. 
Conversely, given such a pair $(D, h)$, they give 
charts to $\tilde M$, and hence induces a $(G,X)$-structure on
$\tilde M$. Since a deck-transformation is 
a $(G,X)$-map $\tilde M \ra \tilde M$, we see that $M = \tilde M/\pi_1(M)$
has an induced $(G,X)$-structure from $\tilde M$.

We say that two such pairs $(D, h)$ and $(D', h')$ are 
{\em $G$-equivalent}
if $D' =\vth \circ D$ and $h'(\cdot) = \vth \circ h(\cdot)\circ \vth^{-1}$
for $\vth \in G$. 

Let us look at 
the set $\mathcal{M}(M)$ of all $(G, X)$-structures on $M$ 
and introduce an equivalence relation that two
$(G,X)$-structures $\mu_1$ and $\mu_2$ are equivalent 
if there is an isotopy $\phi:M \ra M$ so that 
the induced $(G, X)$-structure $\phi^*(\mu_1)$ obtained 
by pulling back charts equals $\mu_2$.
The {\em deformation space of $(G,X)$-structures on $M$} 
(without topology) is defined to be this set $\mathcal{M}(M)/\sim$. 

We reinterpret this space as follows:
consider the set of diffeomorphisms $f:M \ra M'$ 
where $M'$ is a $(G,X)$-manifold. We introduce 
an equivalence relation that $f$ and $f':M \ra M''$ 
are equivalent if there is a $(G,X)$-diffeomorphism 
$\phi:M' \ra M''$ so that $\phi\circ f$ is isotopic
to $f'$. The set of equivalence classes corresponds 
in one-to-one manner with the above 
space by sending $f:M \ra M'$ to $f_*(\mu)$ for
the $(G,X)$-structure $\mu$ on $M'$.

We present yet another version of this set:
Clearly, $\tilde M \times I$ is a universal cover 
of $M \times I$ and the group of deck transformation 
group is isomorphic to $\pi_1(M)$ with an obvious action.
We identify $\pi_1(M \times I)$ with $\pi_1(M)$. 
Consider the set of diffeomorphisms 
$\tilde f:\tilde M \ra \tilde M'$ 
equivariant with respect to
an isomorphism $\pi_1(M) \ra \pi_1(M')$
where $M'$ is a $(G,X)$-manifold. 
Introduce an equivalence relation that 
$\tilde f:\tilde M \ra \tilde M'$
and $\tilde f':\tilde M \ra \tilde M''$ 
are equivalent if there are a $(G,X)$-diffeomorphism
$\tilde \phi:\tilde M' \ra \tilde M''$ and an isotopy
$H:\tilde M \times [0,1] \ra \tilde M''$ 
equivariant with respect to the isomorphisms
\[\tilde \phi_*:\pi_1(M') \ra \pi_1(M'') \hbox{ and } 
\tilde f'_*: \pi_1(M) \ra \pi_1(M'')\] respectively
so that 
\[H_0 = \phi\circ \tilde f 
\hbox{ and } H_1 = \tilde f'.\]
Let us denote this space by $D_I(M)$.
The set of the equivalence classes is certainly in 
one-to-one correspondence with the above set
since two different choices of lifts of
$f:M \ra M'$ differ by a deck transformation
of $\tilde M'$ which is a $(G,X)$-diffeomorphism.

We now give the final version in order to introduce topology:
Following Lok's thesis \cite{Lok:84}, 
we define the {\em isotopy-equivalence 
space $\mathcal{S}(M_0)$ of $(G,X)$-structures} for 
a compact orbifold $M_0$ 
to be the space of equivalence classes of pairs
$(D, \tilde f:\tilde M_0 \ra \tilde M)$ where $\tilde f$ is a diffeomorphism
equivariant with respect to an isomorphism $\pi_1(M_0) \ra \pi_1(M)$
and $D:\tilde M \ra X$ is an immersion equivariant with respect to
a homomorphism $\pi_1(M) \ra G$.
Two such pairs $(D, \tilde f)$ and $(D', \tilde f')$ 
are {\em isotopy-equivalent} if 
and only if there are a diffeomorphism 
$\tilde \phi:\tilde M \ra \tilde M'$ lifting 
a diffeomorphism $\phi:M \ra M'$
with $D'\circ \tilde \phi = D$
and an isotopy $H: \tilde M\times [0,1] \ra \tilde M'$ 
equivariant with respect to the isomorphism 
$\tilde \phi_*:\pi_1(M) \ra \pi_1(M')$ 
so that $H_0 = \tilde \phi\circ \tilde f$ and $H_1 = \tilde f'$.

The topology is given on the set of pairs 
by $C^s$-topology on $D\circ \tilde f$, i.e., 
a sequence of functions converges if it does on every compact subset
of $\tilde M_0$ uniformly in $C^s$-sense. 
($s \geq 1$ is sufficient for all purposes.)
We give the quotient topology on $\mathcal{S}(M_0)$. 

There is a natural $G$-action on $\mathcal{S}(M_0)$ given
by 
\[\gamma (D, \tilde f) = (\gamma \circ D, \tilde f), \gamma \in G.\]
Let $\mathcal{D}(M_0)$ be the quotient space under this action. 
Then $D_I(M_0)$ and $\mathcal{D}(M_0)$ are also in one-to-one
correspondence given by sending $\tilde f:\tilde M_0 \ra \tilde M'$ 
to the equivalence class of $(D, \tilde f)$ where 
$D:\tilde M' \ra X$ is a developing map of $M'$.
Therefore, we call $\mathcal{D}(M_0)$ the {\em deformation space} 
of $(G,X)$-structures on $M_0$. 

The set of all homomorphisms $h:\pi_1(M) \ra G$ is denoted
by $\Hom(\pi_1(M), G)$. We restrict our attention 
to the case that $\pi_1(M)$ is finitely presented. 
Let $g_1, \dots, g_n$ denote
the generators of $\pi_1(M)$, and $R_1, \dots, R_m$ 
the relations. Then $H=\Hom(\pi_1(M), G)$ can be injectively
mapped into $G^n$ by sending a homomorphism $h$ to the element
$(h(g_1),\cdots, h(g_n))$ corresponding to generators. 
The relations give us
the subset of $G^n$ where $H$ can lie.
Actually, the subset defined by the relations gives us
precisely the image. Thus, we identify $H$
with this subset. The subset has a subspace topology
of a real algebraic set, which we give to $H$. 
Obviously, it is a metric space if $G$ has a metric.

There is an action by conjugations on $H$ sending 
a homomorphism $h(\cdot)$ to $\vth \circ h(\cdot) \circ \vth^{-1}$ for
$\vth \in G$. $H/G$ may not be a Hausdorff space.
There is a subset $H^{s}$ of $H$, where $G$ acts properly,
consisting of points lying in stable orbits when $G$ is the group of
$\bR$-points of an algebraic group $\bar G$ defined over $\bR$.  
$H^{s}/G$ is a Hausdorff real analytic space.

We define a {\em pre-holonomy map} 
\[\mathcal{PH}: \mathcal{S}(M_0) \ra \Hom(\pi_1(M_0), G)\]
by sending $(D, \tilde f:\tilde M_0 \ra \tilde M)$ to the holonomy
representation $h\circ \tilde f_*$ where $\tilde f_*$ is the induced 
homomorphism $\pi_1(M_0) \ra \pi_1(M)$.

First of all, this is well-defined:
Let $(D', \tilde f':\tilde M_0 \ra \tilde M)$ be 
an equivalent pair. Then $D = D' \circ \tilde \phi$ for
a lift $\tilde \phi$ of an isotopy $\phi:M \ra M'$.
For a deck transformation $\vth$ of $\tilde M_0$,
we obtain
\begin{eqnarray} 
h (\tilde f_*(\vth))\circ D  
&=& D \circ \tilde f \circ \vth \circ \tilde f^{-1}
\nonumber \\
&=& D'\circ \tilde \phi \circ \tilde f \circ 
\vth \circ \tilde f^{-1} \circ \tilde \phi^{-1} \circ 
\tilde \phi \nonumber \\ 
&=& D' \circ \tilde f' \circ \vth \circ 
\tilde f^{\prime -1} \circ \tilde \phi 
\quad \quad \mbox{ by Proposition \ref{prop:diffeomorphisms} } \nonumber \\ 
&=& h'(\tilde f'_*(\vth))\circ D'\circ \tilde \phi \nonumber \\
&=& h'(\tilde f'_*(\vth))\circ D.  
\end{eqnarray} 
Therefore, we obtain
$h \circ \tilde f_* = h'\circ \tilde f'_*$. 

Also $\mathcal{PH}$ is continuous: Let $C(M_0)$ denote
the space of pairs $(D, \tilde f:\tilde M_0\ra \tilde M)$ with 
$C^s$-topology.
Then sending $(D, \tilde f:M_0 \ra M)$ to $h\circ \tilde f_*$ is
a continuous map, which we denote by 
\[\mathcal{PPH}:C(M_0) \ra \Hom(\pi_1(M), G)\]
for later purposes:
the $C^s$-convergence of sequence of $D_i\circ \tilde f_i$
restricted on compact subsets 
for a sequence of pairs $(D_i, \tilde f_i)$ implies the uniform 
$C^\infty$-convergence of $h_i(\tilde f_*(\vth))$ for 
each deck transformation $\vth$. 
(The sequence of locally defined maps 
$D_i\circ \tilde f_i \circ \vth \circ \tilde f_i^{-1} \circ D_i^{-1}$  
converges in $C^s$-topology 
in sufficiently small compact domains in $X$
and hence in $C^\infty$-topology
as $G$ acts smoothly on $X$. $M_0$ needs to be compact here.)

\section{The proof of Theorem 1}

Again, let $G$ be a Lie group acting on a space $X$ 
smoothly with the local properties mentioned above.
Let us now present three lemmas \ref{lem:pertu},
\ref{lem:pertu2}, and \ref{lem:Gpertu}
on the perturbation of the finite group actions and 
conjugation by diffeomorphisms.
\begin{lem}\label{lem:pertu}
Let $G_B$ be a finite subgroup of $G$ acting
on an $n$-ball $B$ in $X$. Let $h_t:G_B \ra G$, 
$t\in [0, \eps]$, $\eps >0$, be an analytic 
parameter of representations of $G_B$ so
that $h_0$ is the inclusion map.
Then for $0 \leq t \leq \eps$, 
there exists a continuous family of
diffeomorphisms $f_t: B \ra B_t$ 
to open balls $B_t$ in $X$ so that $f_t$ conjugates 
$h(G_B)$-action to $h_t(G_B)$-action\/{\rm ;} i.e., 
$f_t^{-1}h_t(g)f_t = h(g)$ for each $g \in G_B$
and $t \in [0, \eps]$.
\end{lem}
\begin{proof} 
We take a product $X \times [0,1]$ and let 
$v$ be a vector field in the positive $[0,1]$-direction 
in the product space. Then $G_B$ acts smoothly on
$X \times [0,1]$ by sending $(x, t)$ to $(h_t(g)(x), t)$
for $g \in G_B$. 
We average $g^*(v)$ for $g \in G_B$ to obtain
a smooth $G_B$-invariant vector field $V$. The integral
curve $l$ of $V$ starting from $(x, 0)$ is mapped
to an integral curve $m$ of $V$ starting from
$(g(x), 0)$ by the $G_B$-action. 
Thus, the endpoint $l(1)$ is sent
to $m(1)$, and so $g(l(1)) = m(1)$. 
Hence, let $f'_t(x)$ equal the point of
the path from $(x,0)$ at time $t$.
Then $f_t = p\circ f'_t:X \ra X$ for the projection
$p:X \times [0,1] \ra X$ is a desired diffeomorphism
and $f_t(B)$ is the desired open ball.
\end{proof}

A point $x$ of a real algebraic set has a neighborhood
with a semi-algebraic homeomorphism to 
a cone over a semi-algebraic set $S$ in
the boundary of a small ball with a cone-point
at the origin corresponding to $x$.

\begin{lem}\label{lem:pertu2}
Let $G_B$ be a finite subgroup of $G$ acting on 
an $n$-ball $B$. 
Suppose that $h$ is a point of 
an algebraic set $V=\Hom(G_B, G)$ for
a finite group, and let $C$ be a cone-neighborhood
of $h$. Then for each $h'\in C$, there is 
a corresponding diffeomorphism 
\[f_{h'}:B \ra B_{h'}, B_{h'}= f_{h'}(B)\]
so that $f_{h'}$ conjugates 
the $h(G_B)$-action on $B$ to 
the $h'(G_B)$-action on $B_{h'}$\, {\rm ;} 
i.e., $f_{h'}^{-1}h'(g)f_{h'} = h(g)$ for each $g \in G$.
Moreover, the map $h' \mapsto f_{h'}$ is 
continuous from $C$ to the space $C^\infty(B, X)$
of smooth functions from $B$ to $X$. 
\end{lem}
\begin{proof} 
Parameterize $C$ by $[0,\eps]\times S$ for a semi-algebraic
set $S$ with $\{0\}\times S$ 
corresponding to $h$ and, for each $x \in S$, there is 
a map $[0,\eps]\times x \ra C^\infty(B, X)$
from the above lemma \ref{lem:pertu}.
Again, we obtain a smooth $G_B$-invariant vector field
$V_x$ on $X \times [0,\eps]$ as above, 
and $V_x$ depends continuously on $x$. 
From this, we see that $f_{x, t}$ corresponding 
to a representation corresponding to $(x, t)$
depends continuously on $(x, t)$.
\end{proof}

An isotopy of an embedded submanifold extends 
to one of the ambient manifold in a continuous manner,
which is the following version of 
Cerf's ``first isotopy and extension 
theorem'' \cite{Cerf} (see Lok \cite{Lok:84}): 

\begin{thm}\label{thm:cerf} 
Let $F$ be a closed smooth submanifold of $X$ with corners. 
Let ${\mathcal{E}}(X)$ denote the space of 
isotopies $X\times [0,1] \ra X$ with the $C^s$-topology. 
Let ${\mathcal{E}}(F, X)$ denote the space 
of imbeddings of $F$ in $X$ with $C^s$-topology. 
Consider the map 
\[\Phi: {\mathcal{E}}(X) \ra {\mathcal{E}}(F,X)\]
given by sending an isotopy $f_t$ to $f_1|F$. 
Then there is a neighborhood of the inclusion
$i:F \ra X$ of ${\mathcal{E}}(F, X)$ on 
which there is a continuous section 
$s$ of $\Phi$ and $s(i) = e$ 
where $e$ is the identity isotopy
in ${\mathcal{E}}(X)$. 
\end{thm}

Continuing to use the notation of Lemma \ref{lem:pertu2},
we define a parameterization $l: S\times [0,\eps] \ra C$
which is injective except at $S\times \{0\}$ which
maps to $h$. (We fix $l$ although $C$ may become smaller 
and smaller). For $h'\in S$,
we denote by $l(h'):[0,\eps] \ra C$ be a ray in $C$ 
so that $l(h')(0) = h$ and $l(h')(\eps) = h'$.
Let the finite group $G_B$ act on a submanifold $F$ of $B$.
A {\em $G_B$-equivariant isotopy} $H: F\times [0,\eps'] \ra X$
is a map so that $H_t$ is an imbedding 
for each $t \in [0,\eps']$, with $0 < \eps' \leq \eps$, 
conjugating the $G_B$-action on $F$
to the $l(h')(t)(G_B)$-action on $X$, where
$H_0$ is an inclusion map $F \ra X$.
The above lemma \ref{lem:pertu2} says that for each $h'\in C$,
there exists a $G_B$-equivariant isotopy 
$H:B \times [0,\eps] \ra X$. 
We will denote by $H(h')_{\eps'}:B \ra X$ 
the map obtained from $H$ for $h'$ and $t =\eps'$.
Note also by the similar proof,
for each $h'\in S$, there exists 
a $G_B$-equivariant isotopy $H:F \times [0,\eps] \ra X$. 

\begin{lem}\label{lem:Gpertu}
Let $H: F \times [0,\eps'] \times S \ra X$ be a map so
that $H(h'):F \times [0, \eps'] \ra X$ is 
a $G_B$-equivariant isotopy of $F$ for each $h'\in S$
where $0 < \eps' \leq \eps$ for some $\eps > 0$.
Then $H$ can be extended to 
$\hat H: B \times [0, \eps''] \times S \ra X$ 
so that $\hat H(h'):B \times [0, \eps''] \ra X$ 
is a $G_B$-equivariant isotopy of $B$ for each $h'\in S$
where $0 < \eps'' \leq \eps'$.
\end{lem} 
\begin{proof} 
Let ${\mathcal{E}}(C, F, X)$ denote the space of 
all $G_B$-equivariant isotopies of $F \ra X$ 
with the above parametrization $l$.
The maps $H(h')_{\eps'}:F \ra X$ are in 
${\mathcal{E}}(F, X)$. There exists 
a section $s': W \subset {\mathcal{E}}(C, F, X) \ra {\mathcal{E}}(X)$
where $W$ is a neighborhood of $(h, i|F)$
by Theorem \ref{thm:cerf}.
Hence, there exists $\eps''$, $0< \eps'' \leq \eps$,
so that $H(h')_{\delta} \in W$ for 
$0 < \delta < \eps''$ for all $h'\in S$.  
We define $H'(x, \delta, h')$ for $x \in B$ to be 
$s(H(h')_{\delta})_1(x)$. 
This defines a function 
$H': X \times [0, \eps''] \times S \ra X$
so that $H'(h'): X \times [0, \eps''] \ra X$ is 
an isotopy for each $h'$. 

We modify $H'(h')$ to be $G_B$-equivariant. 
We define 
\[H''(h'): X \times [0, \eps''] \ra X \times [0, \eps'']\] 
to be given by $H''(h')(x, t) = (H'(h')(x), t)$. 
Then $H''(h')(x, t)$ for a given $x$ 
is an integral curve of a vector field $V$ on $X\times [0, \eps'']$
with the component in the $[0, \eps'']$-direction equal to $1$.
Since $H(h')$ is an $G_B$-equivariant isotopy, 
we see that $V$ restricted to the image of
$H''(h')(F\times [0, \eps''])$ is $G_B$-invariant. 
We may now average the image vector fields of $V$ 
under $G_B$, and call $V'$ the resulting $G_B$-invariant vector
field on $X\times [0, \eps'']$. Again $V'$ restricted 
to the image equals $V$ on the image
and the second component equals $1$. 
The integral curves of $V'$ give
us a $G_B$-equivariant isotopy 
$\hat H(h'): X \times [0, \eps''] \ra X$ 
extending $H(h'): F \times [0, \eps''] \ra X$. 
Since the section $s$ is continuous, and 
we do averaging and integration, 
it follows that $\hat H: X \times [0, \eps''] \times S 
\ra X$ is continuous. Now restrict $\hat H$ to 
$B \times [0, \eps''] \times S$.
\end{proof} 

\begin{rem}\label{rem:pertu}
We choose some arbitrary Riemannian metric on a neighborhood
of $B$, and can assume that the images of $f_{h'}$s are all
in this neighborhoods (see below). 
By our construction, given any $\eps > 0$, we can make sure
that the $C^s$-norm of $f_{h'}$, constructed in
above lemmas, minus the inclusion map of $B$ is less than $\eps$ 
in some coordinate systems
if we choose the neighborhoods $C$ sufficiently small near $h$. 
In particular, we can assume that for each $\eps >0$, 
there is a neighborhood $C$ of $h$ so that 
$d(f_{h'}(x), x) \leq \eps$ for $x \in B$ and $h' \in C$ 
where $f_{h'}$ is obtained from above three lemmas. 

If $B$ was strictly convex
with smooth boundary, we see that $f_{h'}(B)$ is
also strictly convex with smooth boundary 
as the boundary convexity is given by a $C^2$-condition. 

We can trivially generalize Lemma \ref{lem:Gpertu}  
so that $B$ could be a union of disjoint balls with some finite groups
acting on each.
\end{rem}

Now, we begin the step (I) of the proof of 
Theorem 1 as stated in the outline in the introduction:

Let $M$ be a compact $(G, X)$-orbifold.
We choose a nice finite cover $U_1, \dots, U_k$ of $M$ 
so that $\clo(U_i) \subset W_i$ and 
$\clo(W_i) \subset V_i$ for nice finite covers 
$W_1, \dots, W_k$ and $V_1, \dots, V_k$. 
(This can be done by change radii by small amounts.
See the proof of Proposition \ref{prop:nicecover}.)
We assume that $U_i, W_i, V_i$ all have 
the open balls as the open subsets in the model pairs.

Give $M$ a Riemannian metric, and 
let $\tilde M$ be the universal cover of $M$. 
Since $\tilde M$ is a manifold, it has an induced
Riemannian metric in the ordinary sense. 
The components of the inverse images of the balls $V_i$ above, 
are strictly convex balls which are images of 
exponential maps. By their strict convexity, 
any two of them meet in a strictly convex ball, 
i.e., in a contractible subset. 

For each $V_i$, choose an arbitrary component 
$L_i$ in $\tilde M$ of its inverse image. 
$L_i$ is homeomorphic to an $n$-ball, and
there exists a finite subgroup $\Gamma_i$ of 
the fundamental group $\Gamma$ of $\tilde M$
acting on $L_i$, and $(L_i, \Gamma_i)$ is 
a model pair for $V_i$. 
We choose $M_i$ and $N_i$ in $L_i$ corresponding
to $U_i$ and $W_i$ respectively. 

Given $i, j$, if $V_i$ and $V_j$ meet,
then there exists a deck-transformation
$\gamma_{ij}$ so that $L_i \cap \gamma_{ij} L_j \ne \emp$. 
The choice of $\gamma_{ij}$ is not unique
if $\Gamma_i$ and $\Gamma_j$ are not trivial since one can always
multiply $\gamma_{ij}$ in the left by an element of $\Gamma_i$
in the right by an element of $\Gamma_j$. 
Let $\Gamma_{ij}$ denote the all such possibilities
for $L_i$ and $L_j$. Clearly, we have
\begin{equation}\label{eqn:gammaji}
\gamma^{-1} \in \Gamma_{ji} \hbox{ if
and only if } \gamma \in \Gamma_{ij}. 
\end{equation}

If $L_i \cap \gamma(L_i) \ne \emp$,
then $\gamma$ acts on $L_i$ since $L_i$ is a normal
neighborhood. Hence, we have 
\begin{equation}\label{eqn:gammaii}
\Gamma_{ii} =\Gamma_i.
\end{equation}

Clearly, every element of $\Gamma_{ij}$
can be written $\gamma_1\gamma\gamma_2$ 
where $\gamma_1 \in \Gamma_i$, $\gamma_2 \in \Gamma_j$, 
and $\gamma$ is a fixed element of $\Gamma_{ij}$.
Thus, one can make sense of 
the statement that the coset space
$\Gamma_{ij}/\Gamma_j$ is in one-to-one correspondence with
$\Gamma_i$. 

We note that $L_i\cap \gamma L_j$ for $\gamma\in \Gamma_{ij}$
is a convex ball, hence contractible. 
The same can be said for $M_i\cap \gamma M_j$
and $N_i \cap \gamma N_j$. 
We assume that $L_i \cap \gamma L_j \ne \emp$ 
if and only if $M_i \cap \gamma M_j \ne \emp$ 
if and only if $N_i \cap \gamma N_j \ne \emp$. 

We claim that $\bigcup_{i,j} \Gamma_{ij}$ is 
a set of generators of $\pi_1(M)$: 
Let $\gamma \in \pi_1(M)$. Since $\tilde M$ is connected,
there is a path from $L_1$ to $\gamma(L_1)$. 
There exists a collection 
$A_1, A_2, \dots, A_n$ of open sets 
so that $A_1 = L_1$, $A_n = \gamma(L_1)$,
$A_j \cap A_{j+1} \ne \emp$ for $j = 1, \dots, n-1$, 
and $A_j$ is of form $\gamma_j(L_{k_j})$ for some $k_j$ and
$j=1, \dots, n$. 
Since $A_j$ and $A_{j+1}$ meet, 
and so $\gamma_j(L_{k_j})$ and $\gamma_{j+1}(L_{k_{j+1}})$ 
meet, it follows that $\gamma_j^{-1}\gamma_{j+1}$ 
lies in $\Gamma_{k_j, k_{j+1}}$. 
We have 
\begin{eqnarray}
A_1 &=& L_1 \nonumber \\
A_2 &=& \gamma_{1{k_2}} L_{k_2}, 
\mbox{ for } \gamma_{1{k_2}} \in \Gamma_{1{k_2}} \nonumber \\
A_3 &=& \gamma_{1{k_2}} \gamma_{k_2 k_3} L_{k_3},
\mbox{ for } \gamma_{k_2 k_3} \in \Gamma_{k_2 k_3} \nonumber \\ 
\vdots & \vdots & \vdots \nonumber \\ 
A_n = \gamma(L_1) &=& \gamma_{1{k_2}} \gamma_{k_2 k_3} \cdots 
\gamma_{k_{n-1}1} L_1, 
\mbox{ for } \gamma_{k_{n-1}1} \in \Gamma_{k_{n-1} 1} 
\end{eqnarray}
Thus, we see that 
\[\gamma = \gamma'\gamma_{1{k_2}} \gamma_{{k_2} {k_3}} \cdots 
\gamma_{k_{n-1}1} \hbox{ for some } \gamma' \in \Gamma_1 \] 
We can write any element of $\Gamma$ as a product of
elements in $\bigcup_{i, j} \Gamma_{ij}$.

Also, we see that 
\begin{equation}\label{eqn:grelation}
\gamma\circ \gamma' \in \Gamma_{ik}
\end{equation}
if $\gamma \in \Gamma_{ij}$ and $\gamma' \in \Gamma_{jk}$
and $\gamma\circ \gamma'(x) \in L_i$ for some 
$x \in L_k$. 

Next, we do the step (II) of the outline. More precisely, 
we will find a neighborhood $\Omega$ of $h\circ \tilde f_*$ 
in $\Hom(\pi_1(M_0),G)$ so that there is 
a continuous map $s:\Omega \ra C(M_0)$ where
$\mathcal{PPH}\circ s$ is the identity map and 
$s(h\circ f_*) = (D, \tilde f)$ where $D$ is a developing map
$\tilde M \ra X$ where $\tilde M$ is a universal cover 
of an orbifold $M$, and $\tilde f$ is a lift of a diffeomorphism 
$M_0 \ra M$. The map $s$ induces a continuous map 
\[\tilde s: \Omega \ra \mathcal{S}(M_0),\]
which is a local section of $\mathcal{PH}$. 

This will be accomplished by the following steps: 
(i) First we specify $\Omega$. Perturbations in $\Omega$ 
induce deformations of the model 
neighborhoods, and we construct the orbifold using 
the deformations in $\Omega$ of the holonomy 
of the patching deck transformations.
(ii) We show that the constructed orbifold 
is diffeomorphic to $M_0$. This will be done 
by patching together the deformation maps of
model neighborhoods. (iii) We show that 
the diffeomorphism lifts to the diffeomorphisms 
of the universal cover. Using this fact, 
we can show that the deformed orbifold indeed
has the desired deformed holonomy homomorphism.

(i)
Let $(D, \tilde f:\tilde M_0 \ra \tilde M)$ be an element of $C(M_0)$.
Let $h$ be the associated holonomy homomorphism
$\pi_1(M_0) \ra G$. 

One can construct the underlying space 
of $X_M$ from $V_i$s.
That is, we introduce an equivalence relation 
on the disjoint union 
$\coprod_{i=1}^n L_i$ given by letting
$x \sim y$ if $x =\gamma_{ij}(y)$ for 
$x \in L_i, y \in L_j$. Obviously, the orbifold
structure is encoded in this construction;
thus, we can construct $M$ back from these. 

We can also construct $M$ from
$\coprod_{i=1}^n D(L_i)$ from the equivalence
relation that $x \sim_M y$ if 
$x \in h(\gamma_{ij})(y)$ for $x \in D(L_i),
y \in D(L_j)$. This is easily shown to be an equivalence
relation (see equations (\ref{eqn:gammaji}), 
(\ref{eqn:gammaii}), and (\ref{eqn:grelation})).
Let $Q: \coprod D(L_i) \ra M$ denote
the quotient map. The components of 
the inverse images of $Q(D(L_i))$ under 
the universal covering map $\tilde M \ra M$ 
form a covering of $\tilde M$. 

\begin{rem}\label{rem:stab}
The open sets of form $\gamma(L_i)$ constitute a cover
of the universal covering orbifold $\tilde M$. For given 
three sets $\gamma_{i_1}(L_{i_1})$,
$\gamma_{i_2}(L_{i_2})$, and $\gamma_{i_3}(L_{i_3})$ 
so that $\gamma_{i_l}(L_{i_l}) \cap \gamma_{i_{l+1}}(L_{i_{l+1}}) 
\ne \emp$ with $l=1,2,3$ in cyclic sense, 
we require that $D$ restricted to 
their union should be an imbedding
and their intersection should be of generic type
in $C^\infty$ deformations of $L_i$s. 
We require the same pattern for $M_i$ and $N_i$ as well.
(We don't want a sudden change in the intersection 
pattern of these three sets, i.e., we need 
the stability.)
\end{rem}

For convenience, we identify $M$ with $M_0$ and 
$\pi_1(M)$ with $\pi_1(M_0)$ by 
$\tilde f:\tilde M_0 \ra \tilde M$:
We choose a cone-neighborhood $\Omega$ of $h$ 
in $\Hom(\pi_1(M), G)$ so that for each finite group
$\Gamma_i$ associated with $L_i$ are in a neighborhood $\hat C_i$ 
of $\Hom(\Gamma_i, G)$ satisfying Remark \ref{rem:pertu}
for $B$ equal to $D(N_i)$ or $D(M_i)$,  
and we choose Riemannian metrics on $D(L_i)$ from $M$ 
pushed to $X$ by the map $D| L_i$. 
(Also, we fix a parameterization
of $\hat C_i$ by $[0,\eps] \times S_i$ for some semi-algebraic 
set $S_i$.)
That is, we assume that for $h'$ in $\hat C_i$s, 
the closures of 
$f_{h'}(B)$ are subsets of $D(L_i)$ or $D(N_i)$ for 
$f_{h'}$ obtained in the lemmas \ref{lem:pertu},
\ref{lem:pertu2}, and \ref{lem:Gpertu} respectively.
(In the following, $\hat C_i$ will be modified 
further in various steps; $\Omega$ will be modified 
correspondingly.)

From now on, we will denote by 
the same symbol $f_{h'}$ these functions
for $D(N_i)$ and $D(M_i)$. Also, we denote by 
$D'$ the maps $f_{h'}\circ D$ restricted on $N_i$ and $M_i$
respectively.

Given $h'$ in $\Omega$, we will construct 
a real projective manifold $M'$ which is homeomorphic
to $M_0$. 

We define $D'$ on sets of form $\gamma(N_i)$ or 
$\gamma(M_i)$ for a deck transformation $\gamma$ 
to be $h'(\gamma)\circ D'\circ \gamma^{-1}$ on these sets. 
(They are not yet consistently defined.) 
We need to choose $\Omega$ sufficiently small
so that for sets $\gamma_{i_1}(N_{i_1}), 
\gamma_{i_2}(N_{i_2}), \gamma_{i_3}(N_{i_3})$ 
so that $\gamma_{i_l}(N_{i_l}) \cap 
\gamma_{i_{l+1}}(N_{i_{l+1}}) \ne \emp$ for $l=1,2,3$ 
in cyclic sense, their intersection
pattern does not change under $D'$ 
(as well as under $D$).
Such $\Omega$ exists by Remark \ref{rem:stab}.

We define a topological space 
\[\coprod_{j\in J} D'(N_j)/\sim_{M'}\]
where $\sim_{M'}$ is defined as follows:
$x \in D'(N_i)$ and $y \in D'(N_j)$ are equivalent 
if $x = h'(\gamma)(y)$ for $\gamma \in \Gamma_{ij}$. 
This obviously is reflexive, symmetric, and transitive 
by equations (\ref{eqn:gammaji}), (\ref{eqn:gammaii}),
and (\ref{eqn:grelation}) and the stability.
Let $Q': \coprod D'(N_i)/\sim_{M'} \ra M'$ be the quotient map.

We claim that $M'$ is an orbifold:
We show that $M'$ is Hausdorff. Let $x \in D'(N_i)$ and 
$y \in D'(N_j)$, and suppose that they are not equivalent.
If $i\ne j$ and $\Gamma_{ij} = \emp$, then 
$Q'(D'(N_i))$ and $Q'(D'(N_j))$ are disjoint neighborhoods 
of $Q'(x)$ and $Q'(y)$ respectively.
If $i\ne j$ and $\Gamma_{ij} \ne \emp$, then 
define a map 
\[D'': D'(N_i) \coprod
\coprod_{[\gamma] \in \Gamma_{ij}/\Gamma_j} 
D'(N_j)^{[\gamma]} \ra X\]
by letting $D''|D'(N_i)$ be the inclusion map, 
and $D''|D'(N_j)^{[\gamma]}$ be the map 
$h'(\gamma)|D'(N_j)$ where $\gamma$ is a representative 
of $[\gamma]$ and $D'(N_j)^{[\gamma]}$ 
is a copy of $D'(N_j)$ for each $[\gamma] \in \Gamma_{ij}/ \Gamma_j$. 
(By equation (\ref{eqn:gammaii}), $\Gamma_i = \Gamma_{ij}/\Gamma_j$.) 
Since $x$ and $y$ are not equivalent, 
$D''(x)$ and $D''(\gamma' y^{[\gamma]})$ 
for a copy $y^{[\gamma]}$ of $y$,
every $\gamma' \in \Gamma_j$ and 
$[\gamma] \in \Gamma_{ij}/\Gamma_j$ are not equal. 
We assume without loss of generality that the above 
map $D''$ is an imbedding by choosing our neighborhoods 
sufficiently small. 
Thus, there exist disjoint neighborhoods of $D''(x)$
and the set $\{D''(\gamma' y^{[\gamma]})\} \cap D''(N_i)$ which 
is $\Gamma_i$-invariant.
Then the component of the neighborhood containing $y$ 
and that containing $x$ have no equivalent points 
since every equivalence between $D'(L_i)$ and $D'(L_j)$ 
arises from $\Gamma_{ij}$ (see equation (\ref{eqn:grelation})).
The disjoint neighborhoods 
clearly map to disjoint neighborhoods in $M'$ 
by $Q'$. If $i=j$, a similar argument applies. 
Since we need to consider only finitely many $j$ for
each $i$, the quotient space $M'$ is a Hausdorff space. 

Since $D'(\clo(N_i))$ are compact, and we can 
easily define a map from $\coprod D'(N_i)$ 
to $M'$ by extending the quotient map
$\coprod D'(\clo(N_i)) \ra M'$, we see that $M'$ is
compact. Since $M'$ contains a countable dense subset 
clearly, $M'$ is second countable.

Also, $M'$ is obviously a $(G, X)$-orbifold since
we obtained $M'$ by patching together the 
finite subgroup orbits in open subsets of
$X$: $Q'(D'(N_j))$ form an open
cover of $M'$ modeled on the pairs 
$(D'(N_j), \Gamma_j)$.

(ii) We will construct an orbifold-diffeomorphism $\phi:M \ra M'$: 
Define an imbedding $I_i: Q(D(\clo(M_i))) \ra Q'(D'(N_i))$ by 
\[I_i = Q' \circ f_{h'|\Gamma_i} \circ (Q|D((\clo(N_i)))^{-1}\] 
obtained by Lemma \ref{lem:pertu2}
if $\Gamma_i$ is not trivial,
or \[I_i = Q'\circ (Q|D(\clo(N_i)))^{-1}|D(\clo(M_i))\] if
$\Gamma_i$ is trivial. 
(We define $\tilde I_i: D(\clo(M_i)) \ra D'(N_i)$ 
to be $f_{h'|\Gamma_i}$, which covers
the above map.)
The problem is that $I_i$s are not consistently defined 
over the overlaps of $Q(D(\clo(M_i)))$
and hence, we need to modify the map. 
We have an ordering $M_1, M_2, \dots, M_n$ for
some $n$. We look at the sets of form 
\[Q(D(\clo(M_{i_1})))\cap Q(D(\clo(M_{i_2})))\cap \dots 
\cap Q(D(\clo(M_{i_t}))),\]
with indices satisfying $i_1 < i_2 < \cdots < i_t$
for some $t$. There is an upper bound $t_0$ on $t$. 
(Note that for given $t_0$, the collection 
of the sets of above forms
is composed of quotients of 
disjoint contractible compact submanifolds since our 
covering is nice.)
We define a map $\phi:M\ra M'$ by
defining it to be $I_{i_1}$ on each set of
the above form for $t = t_0$ and the lowest index $i_1$. 
Note that $\tilde I_{i_1}$ defined on
the inverse image of the above set in $D(\clo(M_{i_1}))$ 
is a $\Gamma_{i_1}$-equivariant isotopy. 
Also, since $\phi$ is well-defined,
$\phi$ lifts to a $\Gamma_{i_j}$-equivariant isotopy 
defined on the inverse image of the set in 
$D(\clo(M_{i_j}))$ for each $j=2,\dots, t$ 
mapping to $D'(N_{i_j})$.

We begin an inductive definition:
Suppose that we defined an immersion $\phi$ from 
the union of sets of form 
\begin{eqnarray}
& &Q(D(\clo(M_{i_1})))\cap Q(D(\clo(M_{i_2})))\cap \dots 
\cap Q(D(\clo(M_{i_t}))) \nonumber\\
& & \quad \mbox{ to } Q'(D'(N_{i_1}))\cap Q'(D'(N_{i_2}))\cap \dots 
\cap Q'(D'(N_{i_t})) \subset M'
\end{eqnarray}
for indices satisfying $i_1 < i_2 < \cdots < i_t$
so that $\phi$ lifts to a smooth $\Gamma_{i_j}$-equivariant 
isotopy on the inverse image under $Q$
in $D(\clo(M_{i_j}))$ to $D'(N_{i_j})$, $j=1,\dots, t$.

Then we define a map 
from the union of sets of form 
\begin{eqnarray}
& & Q(D(\clo(M_{i_1})))\cap Q(D(\clo(M_{i_2})))\cap \dots 
\cap Q(D(\clo(M_{i_{t-1}})))
 \nonumber \\
&\quad \mbox{ to } & Q'(D'(N_{i_1}))\cap Q'(D'(N_{i_2}))\cap \dots 
\cap Q'(D'(M_{i_{t-1}})) \subset M' 
\end{eqnarray}
for indices satisfying $i_1 < i_2 < \cdots < i_{t-1}$.
Take one of them say $A$
of form 
\[Q(\clo(D(M_{i_1})))\cap Q(D(\clo(M_{i_2})))\cap \dots 
\cap Q(D(\clo(M_{i_{t-1}})))\]
with indices satisfying
$i_1 < i_2 < \cdots < i_{t-1}$.
The subset $\tilde A = Q^{-1}(A) \cap D(\clo(M_{i_1}))$ 
is an imbedded submanifold on which $\Gamma_{i_1}$ acts.
Let $A'$ be the subset 
\begin{eqnarray*}
&\bigcup_{i_t =1}^n \clo(Q(D(M_{i_1})))\cap  &\clo(Q(D(M_{i_2})))\cap \dots \\ 
& & \cap \clo(Q(D(M_{i_{t-1}}))) \cap \clo(Q(D(M_{i_t}))),
\end{eqnarray*}
of $A$ where $i_1 < i_2 < \cdots < i_{t-1}$, $i_t\ne 
i_1, \dots, i_{t-1}$,
and $\phi$ is already defined with above properties on $A'$.
The subset $\tilde A'=Q^{-1}(A') \cap D(\clo(M_{i_1}))$ 
is an imbedded submanifold of $\tilde A$ on which $\Gamma_{i_1}$ 
acts. $\phi$ lifts to $\tilde \phi$ on $\tilde A'$ 
and using Lemma \ref{lem:Gpertu}, we obtain 
a $\Gamma_{i_1}$-equivariant isotopy 
$\tilde \phi:\tilde A \ra D'(N_{i_1})$. This induces 
an imbedding $\phi:A \ra M'$ extended from $A'$. 
Since $\phi$ is well-defined, $\phi$ lifts 
to a $\Gamma_{i_j}$-equivariant isotopy 
from the inverse image
of $A$ in $D(\clo(M_{i_j}))$ for $j = 1, \dots, t-1$ 
to $D'(N_{i_j})$.
(Note here that the neighborhoods $\hat C_i$ are taken
to be smaller and smaller 
because of Lemma \ref{lem:Gpertu} 
in this induction process.
Also, an ambiguity of choice of the lift is taken 
care of by the fact that $\tilde \phi$ should continuously
deform to an identity map; i.e., $\tilde \phi$ is an isotopy.)
Therefore, the map $\phi$
on $A'$ extends to a smooth map $\phi':A \ra M'$. 
We can do this for sets of form $A$ consistently
since they overlap in sets of form $A'$ where 
$\phi$ is already defined. By induction,
we obtain a map $\phi:M\ra M'$. 

Therefore, we defined for each $M_i$
a map $\tilde \phi_i: D(M_i) \ra D'(L_i)$ 
which is a lifting of $\phi$ from 
the model of neighborhoods $Q(D(M_i))$ to that of $Q'(D'(L_i))$.
For $\eta \in \Gamma_{ij}$, 
\begin{equation}\label{eqn:consistant}
h'(\eta)\circ \tilde \phi_j \circ h(\eta^{-1})| D(M_i) \cap \eta D(M_j)
= \tilde \phi_i| D(M_i) \cap \eta D(M_j).
\end{equation} 
essentially by continuity
since $\tilde \phi_j$ descends to a well-defined 
function $\phi_j$ agreeing with $\phi_i$ on $Q(M_i)\cap Q(M_j)$
and the equation holds if $h' = h$ when $\phi_i$ and $\phi_j$
are the inclusion maps.

By construction,
the map $\phi:M \ra M'$ induces 
a smooth map $\phi|M^r: M^r \ra M^{\prime r}$.
By taking a finite open cover of $M$ initially,
so that there are some points which are covered by
the open sets only once, we see that the local degree of
$\phi|M^r: M^r \ra M^{\prime r}$ is equal to one. 
This map is proper and locally diffeomorphic,
and hence, $\phi|M^r: M^r \ra M^{\prime r}$ is a diffeomorphism. 
Therefore, $\phi:M \ra M'$ is an orbifold-diffeomorphism 
(see the proof of Proposition \ref{prop:nonsingular}).

(iii) Since $M$ and $M'$ are orbifold-diffeomorphic, their
universal covers $\tilde M$ and $\tilde M'$ 
are diffeomorphic equivariantly with respect
to an isomorphism $\pi_1(M) \ra \pi_1(M')$. 
We construct $\tilde M'$ explicitly
from $\tilde M$ as follows: $\tilde M$ is 
covered by open sets of form $\gamma M_i$ for
$\gamma \in \pi_1(M)$, $i=1, \dots, n$. The universal cover
$\tilde M$ can be considered as a quotient space of 
\[\coprod_{\gamma \in \pi_1(M)} h(\gamma)(D(M_i))\] 
under the equivalence
relation that 
\[x \in h(\gamma)(D(M_i)) \sim y\in h(\gamma')(D(M_j))\]
if $x = y$ and $\gamma^{-1}\gamma'\in \Gamma_{ij}$
(or $\gamma(M_i)$ and $\gamma'(M_j)$ meet).
Let $\tilde Q: \coprod h(\gamma)(D(M_i)) \ra \tilde M$
denote the quotient map.
(We take distinct copies in the disjoint union 
of $h(\gamma)(D(M_i))$ for
each $\gamma$ unless $\gamma^{-1}\gamma'$ belongs 
to $\Gamma_{i}$, 
in which case, we consider $h(\gamma)(D(M_i))$,
same as $h(\gamma')(D(M_i))$.)

We define $\tilde M'$ as the quotient space of 
$\coprod h'(\gamma) D'(N_i)$ again
with the equivalence relation 
\[x \in h'(\gamma)(D'(N_i)) \sim y \in h'(\gamma')(D'(N_j))\] 
if $x = y$ and $\gamma^{-1}\gamma' \in \Gamma_{ij}$.
(Again, we use the above copying rule.)
Let $\tilde Q': \coprod h'(\gamma)(D'(N_i)) \ra \tilde M'$
denote the quotient map. 
$\tilde M'$ is shown to be a manifold just as $M'$
is shown to be an orbifold.
Since $M$ is good, $\tilde M$ is a simply-connected manifold.
Also, from a nerve consideration, $\tilde M'$ has the same nerve
of covering as $\tilde M$ as sufficiently 
implied by Remark \ref{rem:stab}.
Thus, $\tilde M'$ is a simply-connected manifold. 
We define a map $p_{M'}:\tilde M' \ra M'$
by defining 
\[p_{M'}| \tilde Q'(h'(\gamma)(D'(N_i))): 
\tilde Q'(h'(\gamma)(D'(N_i))) \ra Q'(D'(N_i))\]
by sending a point corresponding 
to $h'(\gamma)(x)$ to $Q'(x)$ for $x \in \coprod D'(N_i)$.
The map $p_{M'}$ is clearly an orbifold-covering map. 
Moreover, $\pi_1(M)$ acts on $\tilde M'$ by 
letting $\vth \in \pi_1(M)$ act by sending
$x \in \tilde Q'(h'(\gamma)D'(N_i))$ 
to a point in $\tilde Q'(h'(\vth)h'(\gamma)(D'(N_i)))$ 
by a map $\tilde Q' \circ h'(\vth)\circ \tilde Q^{\prime -1}$. 
This is a well-defined automorphism of $\tilde M'$,
and $\tilde M'/\pi_1(M)$ is orbifold-diffeomorphic 
to $M'$. (Of course, the covering map 
$p_M:\tilde M \ra M$ and the action of $\pi_1(M)$ on
$\tilde M$ can be defined the same way.)

The above diffeomorphism $\phi$ lifts to 
a diffeomorphism $\tilde \phi:\tilde M \ra \tilde M'$:
We first recall the lift 
$\tilde \phi_i: D(M_i) \ra D'(N_i)$ 
of $\phi: Q(D(M_i)) \ra Q'(D'(N_i))$,
and for $h(\gamma)(D(M_i))$ where $\gamma \in \pi_1(M)$,
we define 
\[\tilde \phi: h(\gamma)(D(M_i)) \ra h'(\gamma)(D'(N_i))\]
by letting $\tilde \phi(x)$ to be 
$h'(\gamma)\circ \tilde \phi_i\circ h(\gamma)^{-1} (x)$. 
This is well-defined: Let $y$ be a point of 
$h(\gamma')(D(M_j))$ for some $j$, $\gamma' \in \pi_1(M)$
so that $x = y$ and $\gamma^{-1}\gamma' \in \Gamma_{ij}$. 
Then 
\[h'(\gamma')\circ \tilde \phi_j \circ h(\gamma')^{-1}(y) 
= h'(\gamma)h'(\gamma^{-1}\gamma')\circ \tilde \phi_j 
\circ h(\gamma^{\prime -1}\gamma)h(\gamma^{-1})(x).\]
By equation (\ref{eqn:consistant}),
the right-hand side of the above equation is
now $h'(\gamma)\circ \tilde \phi_i \circ h(\gamma^{-1})(x)$.
This defines a smooth map $\tilde \phi:\tilde M \ra \tilde M'$,
which is an immersion.

Since $\tilde M$ and $\tilde M'$ have the same nerve of 
open coverings by open balls, we see that $\tilde M'$ is 
a simply connected manifold. Therefore, $\tilde M'$ is 
a universal covering orbifold of $M'$ by Remark \ref{rem:good}.
Since $p_{M'}\circ \tilde \phi = \phi \circ p_{M'}$ clearly,
$\tilde \phi$ is a lift of an orbifold-diffeomorphism 
$\phi$ and hence is an isomorphism by Corollary \ref{cor:univ}.
The above map $\tilde \phi$ is equivariant, i.e.,
\[\tilde \phi \circ \gamma = \gamma \circ \tilde \phi,
\mbox{ for each } \gamma \in \pi_1(M).\] 
Thus, we see that $\tilde M'$ is the universal covering
space of $M'$ and $\pi_1(M)$ and $\pi_1(M')$ are
isomorphic by $\tilde \phi_*$ induced from $\tilde \phi$. 

We define a developing map $D':\tilde M' \ra X$
by defining 
\[D'| \tilde Q'(h'(\gamma)D'(N_i)) =
\tilde Q^{\prime -1}|h'(\gamma)D'(N_i).\] 
This defines a smooth immersion over $\tilde M'$
in a consistent manner. We consider 
$D'\circ \vth$ for $\vth \in \pi_1(M)$. 
Then on $\tilde Q'(D'(N_i))$, it equals 
\[\tilde Q^{\prime -1}\circ 
\tilde Q'\circ h'(\vth)\circ \tilde Q^{\prime -1}\]
which equals
$h'(\vth)\circ \tilde Q^{\prime -1}$.
We obtain $D'\circ \vth = h'(\vth)\circ D'$. 
Therefore, the holonomy homomorphism is
$h':\pi_1(M) \ra G$ under the identification $\pi_1(M') = \pi_1(M)$.
Or equivalently, $h''\circ \tilde \phi_* = h'$ where $h''$ is 
the holonomy homomorphism of $M'$.

To summarize,
for each $h'\in \Omega$, we defined $M'(h')$ 
with a development pair $(D', h'')$ and 
a diffeomorphism $\phi_{h'}:M \ra M'$ lifting to 
a diffeomorphism $\tilde \phi_{h'}: \tilde M \ra \tilde M'(h')$ 
so that $h''\circ \tilde \phi_{h' *} = h'$. 
(For objects we defined above, we attach a suffix $h'$ 
to indicate that they are constructed for $h'$.)
In fact, we constructed a map
$s':\Omega \ra C(M)$ where 
\[s'(h') = (D', \tilde \phi_{h'}:\tilde M \ra \tilde M'(h')).\]
By Lemma \ref{lem:Gpertu}
and our inductive construction, we can verify that 
\[\tilde \phi_{i, h'}: D(M_i) \ra D'(N_i)\]
depends continuously on $h'$, and hence, so does
\[\tilde \phi_{h'}| h(\gamma)(D(M_i)): 
h(\gamma)(D(M_i)) \ra h'(\gamma)(D'(N_i)).\]
This proves the continuity of
section $s$ completing the step (II).

We will show that 
\[\mathcal{PH}:\mathcal{S}(M_0) \ra \Hom(\pi_1(M_0), G)\]
is locally injective; i.e., for each 
$(D, \tilde f:\tilde M_0 \ra \tilde M)$ there is a neighborhood
where $\mathcal{PH}$ is injective. This will 
be the step (III) of the outline.

Again, we identify $\tilde M$ with $\tilde M_0$ by $\tilde f$.
Let us give $M$ a Riemannian metric
with covering by neighborhoods modeled on
$(L_i, \Gamma_i)$, $i=1, \dots, n$, in $\tilde M$ as above.
We choose $M_i$, $N_i$ as above in $L_i$.

Let $\iota_{\tilde M}:\tilde M \ra \tilde M$ denote
the identity map. 
We choose a neighborhood $\mathcal{O}$ of 
$(D, \iota_{\tilde M}:\tilde M \ra \tilde M)$ 
in $C(M)$ so that 
any two elements $(D_1, \tilde f_1:\tilde M \ra \tilde Y_1)$
and $(D_2, \tilde f_2:\tilde M \ra \tilde Y_2)$ 
satisfy that $D_1\circ \tilde f_1$ is 
sufficiently $C^s$-close to $D_2\circ \tilde f_2$ so that 
\begin{eqnarray*}
&D_1\circ \tilde f_1(\clo(M_i)) \subset D_2\circ \tilde f_2(N_i) 
\subset D(L_i) &\hbox{ and } \\
&D_2\circ \tilde f_2(\clo(M_i)) \subset D_1\circ \tilde f_1(N_i)
\subset D(L_i), &  
\end{eqnarray*}
and the corresponding holonomy homomorphisms 
$h_1$ and $h_2$ belong to $\Omega$ for $\Omega$ defined above. 
(We will add two more conditions on $\mathcal{O}$
making it smaller.)

Let $q:C(M) \ra \mathcal{S}(M)$ be 
the quotient map defined above. 
$q$ is an open map since $\mathcal{S}(M)$ is the space
of orbits in $C(M)$ under the action
of the group of isotopies of $\tilde M$.  

We may assume that 
\[{\mathcal{PH}}(q({\mathcal{O}})) 
= {\mathcal{PPH}}({\mathcal{O}}) \subset \Omega \]
by choosing $\mathcal{O}$ sufficiently small.  

We claim that on $q(\mathcal{O})$, which 
is a neighborhood of the equivalence class of 
$(D, \iota_{\tilde M})$ in $\mathcal{S}(M)$,
$\mathcal{PH}|q(\mathcal{O})$ is injective. 

This will prove Theorem 1
since $\mathcal{PH}|q(\mathcal{O})$ has an inverse map $s$
restricted to the image in $\Omega$. 
The image of $\mathcal{PH}|q(\mathcal{O})$
is open since that of $\mathcal{PPH}|\mathcal{O}$
is open. The latter image is open since 
for each point of its image, we can find a small
neighborhood $\Omega'$ in $\Omega$ 
so that a section $s'$ defined
on $\Omega'$ has images in $\mathcal{O}$ as we can
control the $C^s$-norm of conjugating diffeomorphisms
of model sets by the size of the holonomy perturbations
(see Remark \ref{rem:pertu}.)
Thus, $\mathcal{PH}|q(\mathcal{O})$ is a homeomorphism
to an open subset of $\Hom(\pi_1(M), G)$. 

Given $(D_1, \tilde f_1:\tilde M \ra \tilde Y_1)$ and 
$(D_2, \tilde f_2:\tilde M \ra \tilde Y_2)$ 
in $\mathcal{O}$
with the holonomy homomorphisms $h_1\circ \tilde f_1^*$ and 
$h_2\circ \tilde f_2^*$ which are equal, we show that 
$(D_1, \tilde f_1)$ and $(D_2, \tilde f_2)$ are isotopy equivalent. 
We assumed that 
$D_1\circ \tilde f_1(M_i) \subset D_2\circ \tilde f_2(N_i)$
for each $i$.
We start from $\tilde f_1(M_1)$, and lift 
the map $D_1| \tilde f_1(M_1)$ by $D_2^{-1}$ to $\tilde f_2(N_1)$.

We identify $\pi_1(M)$ with $\pi_1(Y_1)$ and $\pi_1(Y_2)$
by $\tilde f_{1*}$ and $\tilde f_{2*}$ respectively.
If $\gamma(M_j)$ meets $M_1$ for $\gamma \in \Gamma_{1j}$, 
then \begin{eqnarray*}
D_1(\tilde f_1 (\gamma(M_j))) & = &
h_1(\gamma)(D_1(\tilde f_1(M_j))) \subset \\
h_2(\gamma)(D_2(\tilde f_2(N_j)))  & = & D_2(\tilde f_2(\gamma(N_j)))
\end{eqnarray*}
since $h_1(\gamma) = h_2(\gamma)$.
We lift $D_1|\tilde f_1(\gamma(M_j))$ by 
$(D_2|\tilde f_2(\gamma(N_j)))^{-1}$ into $\tilde f_2(\gamma(N_j))$. 
By an induction in this manner, we see that we can
lift an immersion $D_1:\tilde Y_1 \ra X$ 
to an immersion $f_{12}:\tilde Y_1 \ra \tilde Y_2$
by $D_2$ so that $D_2 \circ f_{12} = D_1$. 

Since $h_1 = h_2$, considering $\tilde Y_1$ and $\tilde Y_2$
as quotient spaces of the sets of form $h_1(\gamma)D_1(\tilde f_1(M_j))$ 
and $h_2(\gamma)D_2(\tilde f_2(N_j))$, 
this map is also seen to be $\pi_1(M)$-equivariant;
i.e., 
\[f_{12}\circ \tilde f_1 \circ \gamma 
= \gamma \circ f_{12}\circ \tilde f_1, \gamma \in \pi_1(M);\] 
or in other words,
\[f_{12}\circ \tilde f_{1*}(\gamma) = 
\tilde f_{2*}(\gamma) \circ f_{12}
\hbox{ for } \gamma \in \pi_1(M).\]
Thus, 
\begin{equation}\label{eqn:***}
f_{12*}\circ \tilde f_{1*}(\gamma)= \tilde f_{2*}(\gamma).
\end{equation}

We now show that $f_{12}\circ \tilde f_1$ is isotopic
to $\tilde f_2$ by an isotopy $H:\tilde M \times [0,1] \ra \tilde Y_2$
equivariant with respect to 
the homomorphism $f_{2*}:\pi_1(M) \ra \pi_1(Y_2)$.

Let $Y_2$ have the Riemannian metric pushed by
$\tilde f_2$ with distance metric $d_{Y_2}$. 
Then $\tilde f_2$ is an isometry. Since $M$ is compact, 
$f_{12}\circ \tilde f_1:\tilde M \ra \tilde Y_2$ is a map 
so that 
\[d_{Y_2}(\tilde f_2(x), f_{12}\circ \tilde f_1(x)) \leq \eps
\mbox{ for } x \in \clo(M_i)\]
for some small $\eps > 0$. 

We may choose our neighborhood $\mathcal{O}$ in the beginning 
so that $\eps$ may be chosen to be smaller than 
the minimum radius of the normal neighborhoods for every point of 
$M$. Thus, one can find a unique geodesic from 
$\tilde f_2(x)$ to $f_{12}\circ \tilde f_1(x)$
for each $x \in M$. 
For each point $y$ of $\tilde Y_2$, let 
$v$ be an equivalence class of a vector in $T_y \tilde Y_2$ so that 
$\exp_y(v) = f_{12}\circ \tilde f_1\circ \tilde f_2^{-1}(y)$. 
Since $f_{12}\circ \tilde f_1\circ \tilde f_2^{-1}$ 
is a $\pi_1(Y_2)$-equivariant diffeomorphism by equation \ref{eqn:***},
$v$ is a $\pi_1(Y_2)$-invariant vector field.

Let us denote by $E:T(\tilde Y_2) \ra \tilde Y_2 \times \tilde Y_2$
the map given by sending $(z, w)$ to $(z, \exp_y(w))$ 
for $z\in \tilde Y_2$ and $w\in T_z(\tilde Y_2)$. Then $E$ is 
a differentiable map invertible near the diagonal
$\triangle$ in $\tilde Y_2 \times \tilde Y_2$. Let us call $E^{-1}$
the inverse in a neighborhood of $\triangle$.
Since $E^{-1}$ is a smooth map, 
$v$ is a smooth vector field on $\tilde Y_2$.

If we choose $\mathcal{O}$ 
sufficiently near $(D, \iota_{\tilde M})$,
then $v$ is very small so that 
the map $g_t:\tilde Y_2 \ra \tilde Y_2$
defined by $g_t(x) =\exp_x(tv)$ are immersions
for $t\in [0,1]$: We look at the variation of 
the Jacobian from the Jacobian of the identity map.
Each $g_t$ descends to an immersions $g'_t:Y_2 \ra Y_2$ with 
local degree $1$ isotopic to the identity map.
Since $g'_t$ is a proper-map and of local-degree $1$, 
$g'_t$ restricts to a diffeomorphism $Y_2^r \ra Y_2^r$. 
As we showed above, $g'_t$ are orbifold-diffeomorphisms 
by Proposition \ref{prop:nonsingular}
and so are $g_t$ by Corollary \ref{cor:univ}.
Thus, we require this to hold for $\mathcal{O}$.

Let us denote by $H(y, t)$ the point $\exp_y(tv)$ for
$t\in [0,1]$. Then $H$ is a smooth function 
$\tilde Y_2 \times [0,1] \ra \tilde Y_2$ so that 
$H(y, 0) = y$ for $y\in \tilde Y_2$ 
and $H(f_2(x), 1) = f_{12}\circ \tilde f_1(x)$
for every $x \in \tilde Y_1$.
Therefore, $H$ is $\pi_1(Y_2)$-equivariant isotopy between 
$\tilde f_2$ and $f_{12}\circ \tilde f_1$.
since so are $v$ and $v_t$. 
Thus, $\tilde f_2$ and 
$f_{12}\tilde f_1$ are isotopic, 
and $(D_1, \tilde f_1)$ and $(D_2, \tilde f_2)$
are isotopy-equivalent.
\qed

\bibliographystyle{plain}


\end{document}